\RequirePackage[2020-02-02]{latexrelease}

\documentclass[amsmath,amssymb,superscriptaddress,nobibnotes,aps,noshowpacs,centertags,showkeys]{revtex4}

\usepackage[T1]{fontenc}
\usepackage{algorithm}
\usepackage{algpseudocode}
\usepackage{subfigure} 
\usepackage{color}
\usepackage[english]{babel}

\usepackage{amsthm}
\newtheorem{remark}{Remark}
\newtheorem{theorem}{Theorem}

\makeatletter

\providecommand{\LyX}{L\kern-.1667em\lower.25em\hbox{Y}\kern-.125emX\@}
\DeclareRobustCommand*{\lyxarrow}{%
	\@ifstar
	{\leavevmode\,$\triangleleft$\,\allowbreak}
	{\leavevmode\,$\triangleright$\,\allowbreak}}

\@ifundefined{textcolor}{}
{%
	\definecolor{BLACK}{gray}{0}
	\definecolor{WHITE}{gray}{1}
	\definecolor{RED}{rgb}{1,0,0}
	\definecolor{GREEN}{rgb}{0,1,0}
	\definecolor{BLUE}{rgb}{0,0,1}
	\definecolor{CYAN}{cmyk}{1,0,0,0}
	\definecolor{MAGENTA}{cmyk}{0,1,0,0}
	\definecolor{YELLOW}{cmyk}{0,0,1,0}
}

\makeatother

\begin{document}

\title{Numerical Verification of the Convexification Method for A Frequency-Dependent Inverse Scattering Problem with Experimental Data}
\author{\firstname{Thuy}~\surname{Le}}
\email{tle55@uncc.edu}
\affiliation{%
Department of Mathematics and Statistics, University of North Carolina at Charlotte, Charlotte, NC 28223, USA
}%
\author{\firstname{Vo}~\surname{Anh Khoa}}
\email{anhkhoa.vo@famu.edu, vakhoa.hcmus@gmail.com}
\affiliation{%
Department of Mathematics, Florida A\&M University, Tallahassee, FL 32307, USA
}
\author{\firstname{Michael Victor}~\surname{Klibanov}}
\email{mklibanv@uncc.edu}
\affiliation{Department of Mathematics and Statistics, University of North Carolina at Charlotte, Charlotte, NC 28223, USA
}%
\author{\firstname{Loc Hoang}~\surname{Nguyen}}
\email{loc.nguyen@uncc.edu}
\affiliation{%
Department of Mathematics and Statistics, University of North Carolina at Charlotte, Charlotte, NC 28223, USA
}%
\author{\firstname{Grant}~\surname{Bidney}}
\email{gbidney@uncc.edu}
\affiliation{%
	Department of Physics and Optical Science, University of North Carolina at Charlotte, Charlotte, NC 28223, USA
}%
\author{\firstname{Vasily}~\surname{Astratov}}
\email{astratov@uncc.edu}
\affiliation{%
	Department of Physics and Optical Science, University of North Carolina at Charlotte, Charlotte, NC 28223, USA
}%

\begin{abstract}
The reconstruction of physical properties of a medium from boundary measurements, known as inverse scattering problems, presents significant challenges. The present study aims to validate a newly developed convexification method for a 3D coefficient inverse problem in the case of buried unknown objects in a sandbox, using experimental data collected by a microwave scattering facility at The University of North Carolina at Charlotte. Our study considers the formulation of a coupled quasilinear elliptic system based on multiple frequencies. The system can be solved by minimizing a weighted Tikhonov-like functional, which forms our convexification method. Theoretical results related to the convexification are also revisited in this work

\end{abstract}


\maketitle

\section{Introduction}

In this paper, we build upon our prior research and expand on the performance evaluation of our recently developed \emph{globally convergent} \emph{convexification} numerical method for solving a Coefficient Inverse Problem (CIP) for the 3D Helmholtz equation using multiple frequencies. Our research aims to reconstruct the physical characteristics of explosive-like objects that are buried underground, including antipersonnel land mines and improvised explosive devices (IEDs). Thus, our focus is on three key properties: dielectric constants, locations, and the shapes of front surfaces.

One common approach for numerically solving a CIP is to minimize a conventional least squares cost functional, as described in previous literature such as \cite{Chavent,Gonch1,Gonch2}. However, this method has a major drawback - the cost functional is non-convex and often suffers from the issue of multiple local minima and ravines. As a result, gradient-like methods are limited by getting stuck in any local minimum, and any convergence achieved is only guaranteed if the starting point is in close proximity to the correct solution. Therefore, conventional numerical methods for CIPs are generally limited to local convergence.

\textbf{Definition. }\emph{A numerical method for a CIP is referred to as globally convergent if there exists a theorem that guarantees the method will converge to at least one point within a sufficiently small neighborhood of the correct solution without requiring any prior knowledge of the neighborhood.}

The convexification method is globally convergent, meaning that it is guaranteed to produce at least one solution within a sufficiently small neighborhood of the correct solution, without any prior knowledge of that neighborhood. This method is particularly well-suited for the most challenging cases of solving CIPs, whose data are both backscattering and non-overdetermined. In this context, data are considered non-overdetermined if the number $m$ of free variables in the data is equal to the number $n$ of free variables in the unknown coefficient. In this paper, we consider the case where $m=n=3$. It is worth noting that we are not aware of any other numerical methods for solving CIPs with non-overdetermined data at $m=n\geq 2$ that are both based on the minimization of a conventional least squares cost functional and globally convergent according to the definition given above.

The convexification method has proven effective in solving a 3D CIP with a fixed frequency and a point source moving along an interval of a straight line, as demonstrated by both computationally simulated \cite{KhoaKlibanovLoc:SIAMImaging2020} and experimental data \cite{VoKlibanovNguyen:IP2020,Khoaelal:IPSE2020,Klibanov2021}. In this scenario, we were able to accurately determine the first two key criteria: the dielectric constants and locations of the experimental targets. However, imaging the shapes of the targets' front surfaces requires further improvements. For instance, when dealing with more complicated objects, as shown in Figures \ref{fig test 3}, \ref{fig test 4}, and \ref{fig test 5}, the previous configuration manifests several defects in the reconstructed images. Henceforth, the present paper is focused on further enhancing this aspect.

To address the limitations in imaging the targets' front surfaces with the existing method, we propose to use multiple frequencies while maintaining a fixed point source for the CIP under consideration. This configuration has been previously studied in \cite{LeNguyen:3DArxiv} in conjunction with the convexification technique to solve the same CIP using simulated data. However, its effectiveness with experimental data has only been demonstrated in producing good shapes of objects, while the reconstruction of the dielectric constant is not good. Therefore, we have no choice but to combine this configuration with the previous configuration, which uses a fixed frequency and moving point sources, to amend the third property. In other words,
we have figured out that the best would be to use a two-step procedure. Steps 1 and 2
are performed using two different versions of the convexification method. The version for
Step 1 is described in this paper and the version for Step 3 was described in \cite{KhoaKlibanovLoc:SIAMImaging2020,VoKlibanovNguyen:IP2020,Khoaelal:IPSE2020}.

\textbf{Step 1}. Use the backscattering data for a single location of the
source and multiple frequencies. This gives us accurate geometrical
characteristics of unknown targets: their locations and shapes of front
surfaces. Especially complicated non-convex shapes with voids are imaged
well, see images of letters U, A, O in Figures \ref{fig test 3}--\ref{fig test 5} below. However, values of dielectric
constants of targets are not computed accurately on this step.

\textbf{Step 2}. Use the backscattering data for multiple locations of the
source at a single frequency, as it was done in our previous papers \cite{VoKlibanovNguyen:IP2020,Khoaelal:IPSE2020}.
This provides us with accurate locations and accurate values of dielectric
targets of targets, although the shapes of their front surfaces are not
computed as accurately as they are in Step 1.

\textbf{Step 3}. Assign values of dielectric constants obtained on Step 2 to
images obtained on Step 1. This completes our imaging procedure.

It is worth noting that the configuration of using multiple frequencies and a fixed point source has been studied before in \cite{Liem2018}, but our approach in that study focused on a different approximation procedure using the tail function, rather than the convexification method explored in \cite{LeNguyen:3DArxiv}.

It should be noted that the proposed convexification approach for both the above-mentioned configurations builds upon the ideas of the Bukhgeim--Klibanov method. This method, which is based on Carleman estimates, was initially introduced in 1981 to establish proofs of uniqueness theorems for multidimensional CIPs, as detailed in the seminal work by Bukhgeim and Klibanov \cite{BukhKlib}. Since then, the method has been widely used and extended for solving various inverse problems, see e.g. \cite{Klibanov2013} for a survey of this method. 

The numerical approach considered in this paper deviates from other inversion techniques, such as those employed by Novikov's research group, as described in their publications \cite{Ag, Alekseenko2008, Novikov2015}. These methods address single-frequency data and use distinct treatment methodologies. Additionally, we make reference to \cite{Bakushinskii2020} for a diverse numerical approach to a similar CIP.

The structure of this paper is as follows. In section \ref{sec 2}, we introduce the Coefficient Inverse Problem (CIP) and the corresponding forward problem. Section \ref{sec:3} is devoted to the derivation of our functional $J_{\lambda}$ and the presentation of our theoretical results, which are based on our recent publication \cite{LeNguyen:3DArxiv}. Then, our experimental findings are provided in section \ref{sec num}. Finally, we close the paper by some concluding remarks in section \ref{sec:last}.

\section{Statements of the forward and inverse problem}\label{sec 2}

While the Maxwell's equations are the primary governing equations for the propagation of electromagnetic waves, our paper employs the Helmholtz equation. This approach is supported by numerical demonstrations presented in the appendix of the paper [28], which establish that the Helmholtz equation effectively characterizes the propagation of a specific component of the electric field. Additionally, our successful experimental findings, as reported in our recent publications [2,3], provide further validation for the use of the Helmholtz equation in this context.

Let $\delta$ be the Dirac function. Consider the following time-harmonic Helmholtz wave equation with $\mathbf{x}=\left( x,y,z\right) \in \mathbb{R}^{3}$. 
\begin{equation}\label{eq:helm2}
	\Delta u+ \omega ^{2}\mu \varepsilon ^{\prime }\left( \mathbf{x}%
	\right) 
	u=-\delta \left( \mathbf{x}-\mathbf{x}_{\alpha }\right) \quad \text{in }%
	\mathbb{R}^{3},\;\text{i}=\sqrt{-1}.
\end{equation}%
Physically, $u=u(\mathbf{x})$ can be interpreted as a component of the electric field $E=\left( E_{x},E_{y},E_{z}\right)$ that corresponds to the non-zero component of the incident field. Specifically, in our case, the incident field is characterized by the voltage $E_{y}$. In our experiments, we measure the backscattering signal of this same component. Additionally, $\omega$ represents the angular frequency in rad/m, while $\mu$ and $\varepsilon^{\prime}\left( \mathbf{x}\right)$ denote the permeability (H/m) and permittivity (F/m) of the medium, respectively. The point source $\mathbf{x}_{\alpha}$ is fixed in this study.

We restrict our settings to non-magnetic targets, which means that the materials under consideration have no magnetic properties, and therefore their relative permeability is equal to one. To be more precise, this implies that the ratio of the permeability of the material to the permeability of free space (i.e., vacuum) is unity. Let $\varepsilon _{0}$ represent the vacuum permittivity and let $\mu _{0}$ denote the vacuum permeability. Consider $k = \omega\sqrt{\mu_0 \varepsilon_0}$, equation (\ref{eq:helm2}) can be rewritten as 
\begin{equation}\label{eq:helm3}
	\Delta u+ k^{2}\frac{\mu }{\mu _{0}}\frac{\varepsilon ^{\prime }(\mathbf{x})}{%
		\varepsilon _{0}} u=-\delta \left( \mathbf{x}-\mathbf{x}%
	_{\alpha }\right) \quad \text{in }\mathbb{R}^{3}.
\end{equation}%

We can now express the spatially distributed dielectric constant as $c(\mathbf{x})=\varepsilon ^{\prime }(\mathbf{x})/\varepsilon_0 $. Using this, the conventional Helmholtz equation follows from (\ref{eq:helm3}) and applying the Sommerfeld radiation condition, we get the following system.
\begin{align}
	& \Delta u+k^{2}c(\mathbf{x}) u=-\delta \left( 
	\mathbf{x}-\mathbf{x}_{\alpha }\right) \quad \text{in }\mathbb{R}^{3},
	\label{eq:helm} \\
	& \lim_{r\rightarrow \infty }r\left( \partial _{r}u-\text{i}ku\right)
	=0\quad \text{for }r=\left\vert \mathbf{x}-\mathbf{x}_{\alpha }\right\vert ,%
	\text{i}=\sqrt{-1}.  \label{eq:somm}
\end{align}

Let us now focus on a rectangular prism $\Omega$, defined as $\left( -R,R\right) \times \left( -R,R\right) \times \left( -b,b\right)$ in $\mathbb{R}^{3}$ for $R,b >0$. This prism serves as our computation domain of interest. Besides, we define the lower side of the prism as the near-field measurement site,
\begin{equation*}
	\Gamma :=\left\{ \mathbf{x}:\left\vert x\right\vert ,\left\vert y\right\vert
	<R,z=-b\right\} .
\end{equation*}

In what follows, we make the assumption that the dielectric constant is smooth and meets the following conditions:
\begin{equation}
		\begin{cases}
			c\left( \mathbf{x}\right) \geq 1 & \text{in }\Omega , \\ 
			c\left( \mathbf{x}\right) =1 & \text{in }\mathbb{R}^{3}\backslash \Omega .%
		\end{cases}
	\label{eq:2}
\end{equation}%
The second equation in (\ref{eq:2}) indicates our assumption that the region outside of the domain $\Omega$ is a vacuum. Next, we consider the line of sources that is parallel to the $x$-axis and exists outside of the closure $\overline{\Omega}$. Mathematically, the following line of sources, denoted as $L_{\text{src}}$, is examined:
\begin{equation}
	L_{\text{src}}:=\left\{ \left( \alpha ,0,-d\right) :a_{1}\leq \alpha \leq
	a_{2}\right\},  \label{Lsrc}
\end{equation}%
where $d>b$ and $a_{1}<a_{2}$. Note in this setting that the distance between the line of sources $L_{\text{src}}$ and the $%
xy$-plane is $d$. With this configuration in place, we can now select and describe the fixed point source. The value of $R,d,b,a_1,a_2,\alpha$ will be specified in our experimental results.

\begin{remark}
	To this end, we define the total wave $u$, incident wave $u_i$, and scattered wave $u_s$. It is worth noting that $u=u_{i}+u_{s}$. Besides, the Sommerfeld radiation condition \eqref{eq:somm} is applied to guarantee the existence and uniqueness results for the Helmholtz equation \eqref{eq:helm}; cf. \cite[Chapter 8]{Colton1992}.
\end{remark}

\begin{remark}
	In our configuration of interest, we arrange to measure the data with multiple frequencies. In this regard, the involved waves $u,u_{i},u_{s}$ are dependent of $k$. Henceforth, in the sequel, we write $u=u(\mathbf{x},k),u_{i}=u_{i}(\mathbf{x},k),u_{s} = u_{s}(\mathbf{x},k)$ for $k\in [\underline{k},\overline{k}]$, where $\underline{k},\overline{k}>0$. 
\end{remark}

\subsection{Forward problem}

Prior to introducing the forward problem, we model the incident wave by using the point source,
\begin{equation}
	u_{i}\left( \mathbf{x},k \right) =\frac{\text{exp}\left( \text{i}%
		k\left\vert \mathbf{x}-\mathbf{x}_{\alpha }\right\vert \right) }{4\pi
		\left\vert \mathbf{x}-\mathbf{x}_{\alpha }\right\vert }.  \label{eq:3}
\end{equation}%
We observe that the incident wave $u_i$ satisfies the Helmholtz equation in the form of (\ref{eq:somm}) with $c(\mathbf{x})=1$. By subtracting (\ref{eq:somm}) from the Helmholtz equation for $u_i$, we can obtain the PDE for the scattered wave $u_s$ as follows.
\[
\Delta u_{s}+k^{2}u_{s}=-k^{2}\left(c\left(\mathbf{x}\right)-1\right)u.
\]
Cf. \cite{Colton1992}, the scattered wave is solved via the following integral equation:
\begin{align}\label{eq:4}
		u_{s}\left(\mathbf{x},k\right) & =k^{2}\int_{\mathbb{R}^{3}}\frac{\text{exp}\left(\text{i}k\left\vert \mathbf{x}-\mathbf{x}^{\prime}\right\vert \right)}{4\pi\left\vert \mathbf{x}-\mathbf{x}^{\prime}\right\vert }\left(c\left(\mathbf{x}^{\prime}\right)-1\right)u\left(\mathbf{x}{}^{\prime},k\right)d\mathbf{x}^{\prime}\\
		& =k^{2}\int_{\Omega}\frac{\text{exp}\left(\text{i}k\left\vert \mathbf{x}-\mathbf{x}^{\prime}\right\vert \right)}{4\pi\left\vert \mathbf{x}-\mathbf{x}^{\prime}\right\vert }\left(c\left(\mathbf{x}^{\prime}\right)-1\right)u\left(\mathbf{x}{}^{\prime},k\right)d\mathbf{x}^{\prime},\quad\text{\ensuremath{\mathbf{x}}}\in\mathbb{R}^{3}, \nonumber 
\end{align}%
where we have used the fact that $c-1$ is compactly supported in $\Omega$; see (\ref{eq:2}). Combining \eqref{eq:3} and \eqref{eq:4}, we arrive at the so-called Lippmann-Schwinger equation:
\begin{equation*}
	u\left( \mathbf{x},k \right) =u_{i}\left( \mathbf{x},k \right)
	+k^{2}\int_{\Omega }\frac{\text{exp}\left( \text{i}k\left\vert \mathbf{x}-\mathbf{%
			x}^{\prime }\right\vert \right) }{4\pi \left\vert \mathbf{x}-\mathbf{x}%
		^{\prime }\right\vert } \left( c\left( \mathbf{x}^{\prime
	}\right) -1\right)  u\left( \mathbf{x}^{\prime },k \right) d\mathbf{x}%
	^{\prime },\quad \text{$\mathbf{x}$}\in \mathbb{R}^{3}.
\end{equation*}%

Hence, our forward problem is to determine the boundary information of the total wave field $\left. u\left( \mathbf{x},k \right) \right\vert _{\Gamma }$ for $k\in [\underline{k},\overline{k}]$, based on the known dielectric constant $c$. It is important to remark that the total wave field can be non-zero for all points $x$ in the domain $\Omega$ and for large values of $k$, as demonstrated in \cite{LeNguyen:3DArxiv,KhoaKlibanovLoc:SIAMImaging2020} when $c$ belongs to $C^{15}(\mathbb{R}^3)$ and the Riemannian geodesic line connecting $\mathbf{x}_{\alpha}$ and $\mathbf{x}$ is unique.

\subsection{Coefficient inverse problem (CIP)}

Our CIP is to seek the dielectric constant $c\left( \mathbf{x}\right) $ satisfying (\ref{eq:2}) from knowledge of the boundary measurement $F_{0}\left( \mathbf{x}%
,k\right) $ of the near-field data, 
\begin{equation}
	F_{0}\left( \mathbf{x},k\right) =u\left( \mathbf{x}%
	,k \right) \quad \text{for }\mathbf{x}\in \Gamma ,k\in [\underline{k},\overline{k}],  \label{eq:bdr}
\end{equation}%
where $u\left( \mathbf{x},k \right) $ is the total wave associated with
the incident wave $u_{i}(\mathbf{x},k)$ in \eqref{eq:3}.

While a more detailed description of our experimental setup will be provided in the numerical section, we would like to provide a brief overview. In order to simulate the detection of land mines buried underground, we have buried a single inclusion in a sandbox, with the sand understood as our background medium. The dielectric constant of the sand, $c_{\text{bckgr}}$, is known a priori to be about 4. Although we do not utilize this information in our numerical method, the inclusions in our resulting images are characterized by a dielectric constant $c(\mathbf{x})$ greater than this number 4. It is important to note that the function $c(\mathbf{x})$ used in our mathematical statements incorporates information from both the sand and the inclusion. In order to address this, we measure the raw data twice in our configuration: once when the sandbox is empty, and again when the target is buried within it. By subtracting the former from the latter, we can filter out the information related to the sand. The resulting ``actual'' data (i.e., after subtraction) can then be used in the mathematical setting under consideration.

\begin{figure}[!ht]
	\begin{center}
		\subfigure[Raw data]{\label{far}
			\includegraphics[width=0.4\textwidth]{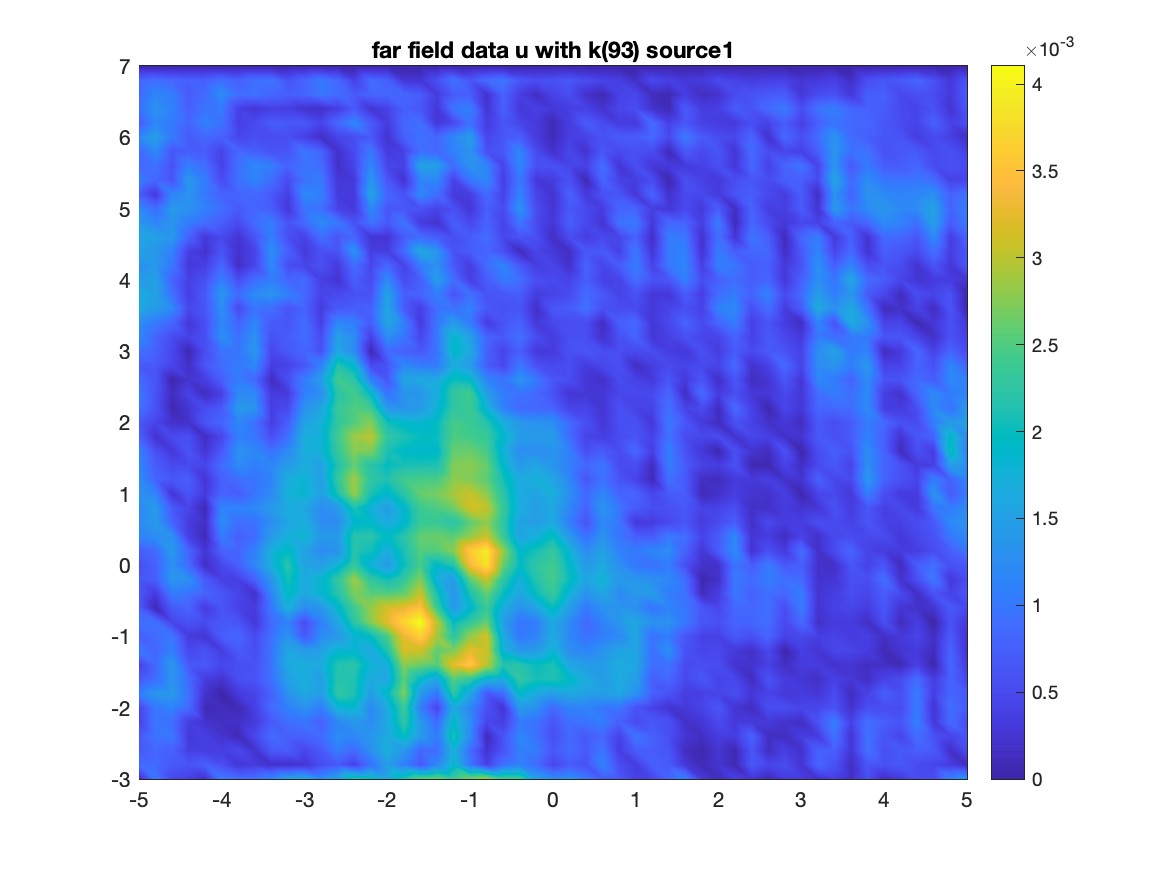} }
		\subfigure[Propagated data]{\label{near}
			\includegraphics[width=0.4\textwidth]{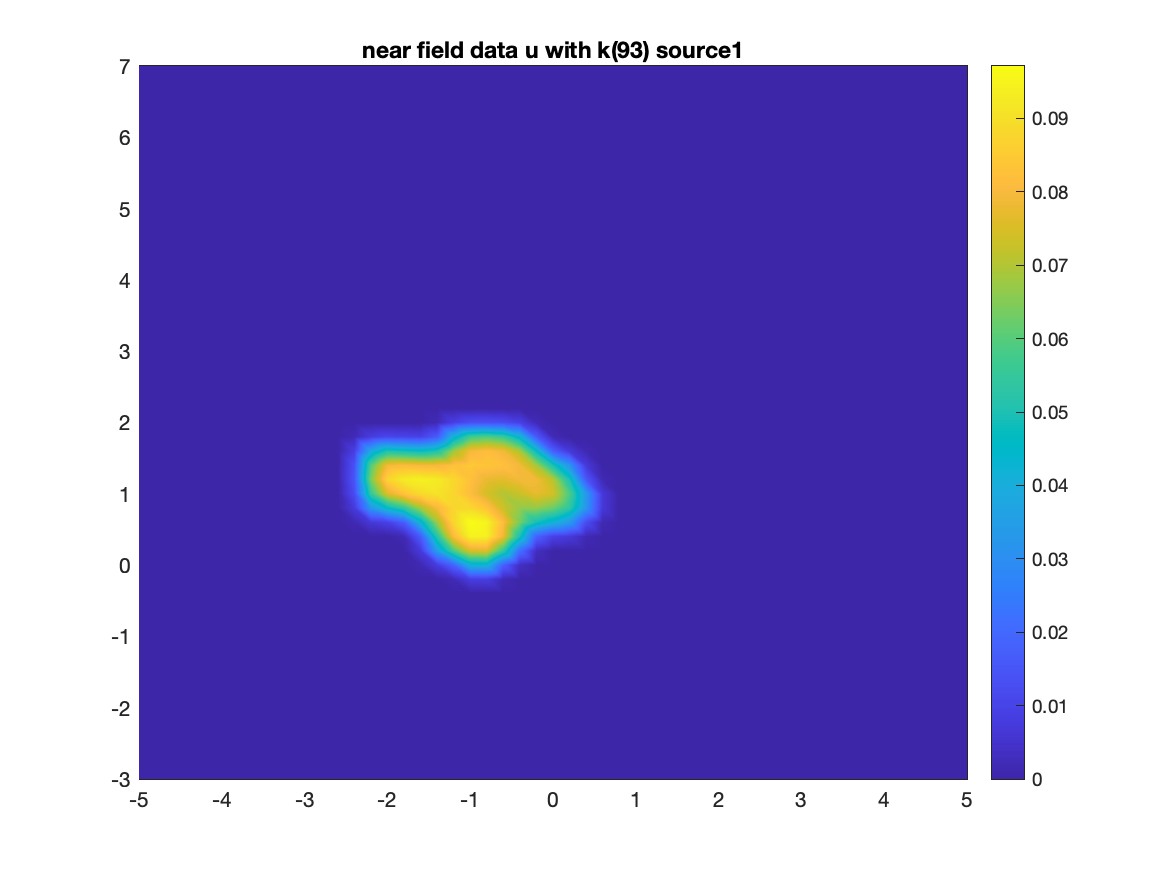}}
	\end{center}
	\caption{Graphical illustration of the absolute value of the raw far-field data (a) and the propagated near-field data (b) in Example 4 in which the experimental object is A-shaped. This set of data is collected at $k = 10.77$ corresponding to $\tilde{f} = 5.14$ (GHz), see \eqref{kless} for the relation between the wavenumber $k$ and the frequency $\tilde{f}$. From these figures, we can see the A shape very clear when using the propagated near-field data. Meanwhile, the shape is not captured well for the raw far-field data.  
	 \label{farnear}}
	
\end{figure}

Our choice to use near-field data stems from our experimental observations, which have shown that far-field data alone do not provide an accurate indication of the buried object's location. In contrast, near-field data have been found to be more reliable, as reported in \cite{Liem2018,VoKlibanovNguyen:IP2020,Khoaelal:IPSE2020} and numerically observed in Figure \ref{farnear}. Furthermore, using near-field data allows us to reduce the size of the computational domain, thereby increasing accuracy, since the number of mesh grids is fixed in our experiments.

Experimentally, we cannot get the near-field data, but the far-field. The ``near-field'' we mean is the approximate dataset that is calculated from the experimental far-field data. To obtain the near-field data, we employ a technique known as data propagation; cf. \cite{Liem2018,VoKlibanovNguyen:IP2020}, the technique is revisited in section \ref{sec num}. This procedure involves eliminating high spatial frequencies to obtain a good approximation of the near-field function $F_0$ in (\ref{eq:bdr}). While we only obtain the measured data $F_0$, our mathematical model requires the $z-$derivative of the function $u\left( \mathbf{x},k \right) $ on $\Gamma$,
\begin{equation}
	F_{1}\left( \mathbf{x},k\right) = \partial_{z}u \left( \mathbf{x}%
	,k \right) \quad \text{for }\mathbf{x}\in \Gamma , k\in [\underline{k},\overline{k}].  \label{1}
\end{equation}%

\begin{remark}
	Given the CIP above, our data are non-overdetermined since the number  $m$ of free variables in the data equals to the number $n$ of free variables in the sought coefficient. In particular, $m = n = 3$ in this scenario.
\end{remark}

\begin{figure}[!ht]
	\begin{center}
		\subfigure[]{\label{conf1}
			\includegraphics[width=0.4\textwidth]{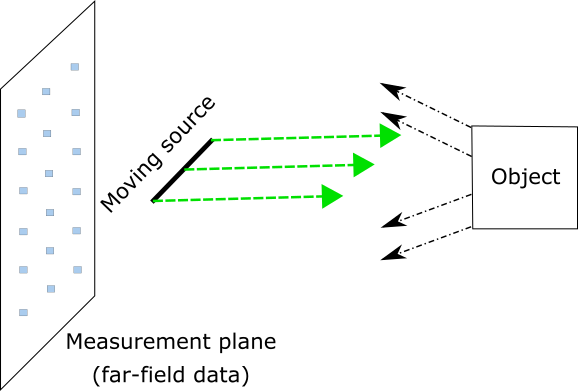} }\hfill
		\subfigure[]{\label{conf2}
			\includegraphics[width=0.4\textwidth]{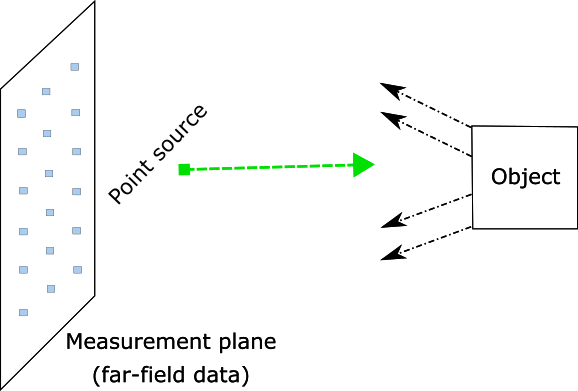}}
	\end{center}
	\caption{(a) Configuration that involves the utilization of multiple sources arranged along a straight line. Each source emits an incident wave at a fixed frequency. (b) Configuration that uses a fixed source while varying the frequency of the generated signal. For both distinct configurations, there are detectors on the measurement plane to collect the backscattering data. \label{diff}}
	
\end{figure}

\begin{remark}
	In the context of multiple sources and a fixed frequency, as presented in e.g. \cite{KhoaKlibanovLoc:SIAMImaging2020}, our CIP can be expressed in a similar way. Specifically, we can consider $\alpha$ as the source variable based on the definition of the line of sources in \eqref{Lsrc}. To handle the configuration of multiple sources and a fixed frequency, we require the boundary measurement $\tilde{F}_{0}(\mathbf{x},\alpha)$, which corresponds to $F_0$ in \eqref{eq:bdr}, and the Neumann-type measurement $\tilde{F}_{1}(\mathbf{x},\alpha)$, which corresponds to $F_1$ in \eqref{1}. Thereby, to derive a coupled quasi-linear elliptic system and establish the convexification, as discussed in section \ref{sec:3}, we can replace the frequency variable $k$ with the source argument $\alpha$ in our formulations. The reader can be referred to Figure \ref{diff} for a visual representation showing the difference in the setup between the two distinct configurations.
\end{remark}

\section{Convexification}\label{sec:3}

\subsection{Derivation of a coupled quasi-linear elliptic system}

For $\mathbf{x}\in \Omega$ and $k\in [\underline{k},\overline{k}]$, we define $
v\left(\mathbf{x},k\right)=\frac{1}{k^{2}}\log\left(\frac{u\left(\mathbf{x},k\right)}{u_{i}\left(\mathbf{x},k\right)}\right)$. It follows from simple calculations that
\begin{align}\label{gradv}
\nabla v\left(\mathbf{x},k\right)=\frac{1}{k^{2}}\left[\frac{\nabla u\left(\mathbf{x},k\right)}{u\left(\mathbf{x},k\right)}-\frac{\nabla u_{i}\left(\mathbf{x},k\right)}{u_{i}\left(\mathbf{x},k\right)}\right],
\end{align}
\begin{align}\label{nablav}
\Delta v\left(\mathbf{x},k\right)=\frac{1}{k^{2}}\left[\frac{\Delta u\left(\mathbf{x},k\right)}{u\left(\mathbf{x},k\right)}-\left(\frac{\nabla u\left(\mathbf{x},k\right)}{u\left(\mathbf{x},k\right)}\right)^{2}-\frac{\Delta u_{i}\left(\mathbf{x},k\right)}{u_{i}\left(\mathbf{x},k\right)}+\left(\frac{\nabla u_{i}\left(\mathbf{x},k\right)}{u_{i}\left(\mathbf{x},k\right)}\right)^{2}\right].
\end{align}
When $\mathbf{x}\in \Omega$, we know that the total wave field $u(\mathbf{x},k)$ satisfies the homogeneous Helmholtz equation, $\Delta u+k^{2}c(\mathbf{x}) u= 0$. The incident wave $u_i(\mathbf{x},k)$ also satisfies this Helmholtz equation with $c(\mathbf{x})=1$. In other words, it holds true that
\[
\frac{\Delta u\left(\mathbf{x},k\right)}{u\left(\mathbf{x},k\right)}=-k^{2}c\left(\mathbf{x}\right),\quad\frac{\Delta u_i\left(\mathbf{x},k\right)}{u_i\left(\mathbf{x},k\right)}=-k^{2}.
\]
Combining this with \eqref{gradv}, \eqref{nablav}, we obtain the following nonlinear PDE for the function $v=v(\mathbf{x},k)$:
\begin{align} \label{PDEv}
\Delta v+k^{2}\left(\nabla v\right)^{2}+\frac{2\nabla v\cdot\nabla u_{i}}{u_{i}}=-c\left(\mathbf{x}\right)+1
\end{align}
for all $\mathbf{x}\in \Omega$ and $k\in [\underline{k},\overline{k}]$. By differentiating \eqref{PDEv} with respect to the argument $k$, we arrive at the following nonlinear PDE:
\begin{align}\label{PDEv2}
\Delta\partial_{k}v+2k^{2}\nabla v\cdot\nabla\partial_{k}v+2k\left(\nabla v\right)^{2}+2\nabla\partial_{k}v\cdot\frac{\nabla u_{i}}{u_{i}}+2\nabla v\partial_{k}\left(\frac{\nabla u_{i}}{u_{i}}\right)=0.
\end{align}

At this stage, it is worth noting that the PDE \eqref{PDEv2} does not contain the unknown function $c(\mathbf{x})$, which is the quantity of interest in our CIP. By solving PDE \eqref{PDEv2}, we can obtain the dielectric constant $c(\mathbf{x})$ via the back-substitution in PDE \eqref{PDEv}.

We, on the other hand, obverse that Equation \eqref{PDEv2} is a non-trivial third-order PDE. Therefore, we rely on the use of a special orthonormal basis of $L^2(\underline{k},\overline{k})$. Denoted by $\left\{\Psi_m\right\}_{m\ge 1}$, this basis is first established in \cite{Klibanov2017}, and it has been applied to solve distinctive inverse problems for PDEs with direct applications to, e.g., electrical impedance tomography and imaging of land mines -- our target application in this work. The reader can be referred to \cite{Klibanov2018,VoKlibanovNguyen:IP2020,Klibanov2023} and references cited therein for an overview of such inverse problems.

To construct this basis, for each $m\ge 1$ we consider  $\varphi_{m}\left(k\right)=k^{m-1}e^{k-(\overline{k}+\underline{k})/2}$. The set
$\left\{ \varphi_{m}\left(k\right)\right\} _{m\ge 1}$ is
linearly independent and complete in $L^{2}\left(\underline{k},\overline{k}\right)$.
We then apply the standard Gram-Schmidt orthonormalization procedure
to obtain the basis $\left\{ \Psi_{m}\left(k\right)\right\} _{m\ge 1}$.

The basis $\left\{ \Psi_{m}\left(k\right)\right\} _{m\ge 1}$ has the following properties:
\begin{itemize}
	\item $\Psi_{m}\in C^{\infty}\left[\underline{k},\overline{k}\right]$ and $\Psi_{m}'$
	is not identically zero for any $m\ge 1$;
	\item Let $S_{mn}=\left\langle \Psi_{n}',\Psi_{m}\right\rangle $ where
	$\left\langle \cdot,\cdot\right\rangle $ denotes the scalar product
	in $L^{2}\left(\underline{k},\overline{k}\right)$. Then the square matrix $\mathcal{S}_{N}=\left(S_{mn}\right)_{m,n=1}^{N}\in\mathbb{R}^{N\times N}$
	is invertible for any $N$ since
\end{itemize}
\[
S_{mn}=\begin{cases}
	1 & \text{if }n=m,\\
	0 & \text{if }n<m.
\end{cases}
\]
Notice that the second property does not hold for either classical orthogonal polynomials or the classical basis of trigonometric functions. The first column of $\mathcal{S}_{N}$ obtained from either of the two conventional bases would be zero.

To solve the third-order nonlinear PDE \eqref{PDEv2}, we consider the truncated Fourier series using the above-mentioned basis. In particular, for $\mathbf{x}\in\Omega$ and $k\in \left[\underline{k},\overline{k}\right]$, we seek
\begin{align}\label{fourier}
v\left(\mathbf{x},k\right)\approx\sum_{n=1}^{N}v_{n}\left(\mathbf{x}\right)\Psi_{n}\left(k\right)=\sum_{n=1}^{N}\int_{\underline{k}}^{\overline{k}}v\left(\mathbf{x},k\right)\Psi_{n}\left(k\right)dk\Psi_{n}\left(k\right).
\end{align}
By plugging \eqref{fourier} into the third-order PDE \eqref{PDEv2}, we find that
\begin{align}
	\sum_{n=1}^{N}\Delta v_{n}\left(\mathbf{x}\right)\Psi_{n}'\left(k\right) & +2\sum_{m,n=1}^{N}\nabla v_{n}\left(\mathbf{x}\right)\cdot\nabla v_{m}\left(\mathbf{x}\right)\left[k^{2}\Psi_{n}\left(k\right)\Psi_{m}'\left(k\right)+k\Psi_{n}\left(k\right)\Psi_{m}\left(k\right)\right] \nonumber\\
	& +2\sum_{n=1}^{N}\nabla v_{n}\left(\mathbf{x}\right)\cdot\left[\Psi_{n}'\left(k\right)\frac{\nabla u_{i}}{u_{i}}+\Psi_{n}\left(k\right)\partial_{k}\left(\frac{\nabla u_{i}}{u_{i}}\right)\right]=0.\label{fourier1}
\end{align}
Henceforth, for $1\le l\le N$ we multiply both sides of \eqref{fourier1} by $\Psi_{l}(k)$ and obtain the following PDE system:
\begin{equation}\label{mainsys1}
	\sum_{n=1}^{N}S_{ln}\Delta v_{n}(\mathbf{x})+\sum_{n,m=1}^{N}P_{lnm}\nabla v_{n}(\mathbf{x})\cdot\nabla v_{m}(\mathbf{x})+\sum_{n=1}^{N}Q_{ln}(\mathbf{x})\cdot\nabla v_{n}(\mathbf{x})=0.
\end{equation}
In (\ref{mainsys1}), we have for $m, n, l = \overline{1,N}$ that
\begin{align*}
	& S_{ln}={\displaystyle \int_{\underline{k}}^{\overline{k}}\Psi_{n}^{\prime}(k)\Psi_{l}(k)dk,}\\
	& P_{lnm}={\displaystyle 2\int_{\underline{k}}^{\overline{k}}\Big(k^{2}\Psi_{n}(k)\Psi_{m}^{\prime}(k)+k\Psi_{n}(k)\Psi_{m}(k)\Big)\Psi_{l}(k)dk,}\\
	& Q_{ln}(\mathbf{x})={\displaystyle 2\int_{\underline{k}}^{\overline{k}}\left(\Psi_{n}^{\prime}(k)\frac{\nabla u_{i}(\mathbf{x},k)}{u_{i}(\mathbf{x},k)}+\Psi_{n}(k)\partial_{k}\frac{\nabla u_{i}(\mathbf{x},k)}{u_{i}(\mathbf{x},k)}\right)\Psi_{l}(k)dk}.
\end{align*}

Now recall from \eqref{eq:bdr} that we measure the wave $u$ on the lower side $\Gamma$ of the prism $\Omega$. Therefore, the Dirichlet boundary information of the sought Fourier coefficients $v_n(\mathbf{x})$ for $1\le n\le N$ is given by
\begin{align}\label{g01}
g_{0n}\left(\mathbf{x}\right)=\int_{\underline{k}}^{\overline{k}}k^{-2}\log\left[F_{0}\left(\mathbf{x},k\right)/u_{i}\left(\mathbf{x},k\right)\right]\Psi_{n}\left(k\right)dk\quad\text{for }\mathbf{x}\in\Gamma.
\end{align}
For $\mathbf{x}\in\partial\Omega\backslash\Gamma$, we apply the heuristic data completion method (cf. e.g. \cite{Liem2018}), choosing that $\left.u\left(\mathbf{x},k\right)\right|_{\Omega\backslash\Gamma}=\left.u_{i}\left(\mathbf{x},k\right)\right|_{\Omega\backslash\Gamma}$. This choice is reasonable because outside of the sandbox is vacuum, i.e. $c=1$. Henceforth, we have
\begin{align}\label{g02}
g_{0n}\left(\mathbf{x}\right)=0\quad\text{for }\mathbf{x}\in\partial\Omega\backslash\Gamma.
\end{align}
As mentioned in the previous section, the Dirichlet measured data \eqref{eq:bdr} can lead to the Neumann-type data \eqref{1}. Moreover, we can compute that for $\mathbf{x}\in \Gamma $,
\[
\partial_{z}v\left(\mathbf{x},k\right)=\frac{1}{k^{2}}\left[\frac{\partial_{z}u\left(\mathbf{x},k\right)}{u\left(\mathbf{x},k\right)}-\frac{\partial_{z}u_{i}\left(\mathbf{x},k\right)}{u_{i}\left(\mathbf{x},k\right)}\right]=\frac{1}{k^{2}}\left[\frac{F_{1}\left(\mathbf{x},k\right)}{F_{0}\left(\mathbf{x},k\right)}-\frac{\partial_{z}u_{i}\left(\mathbf{x},k\right)}{u_{i}\left(\mathbf{x},k\right)}\right].
\]
Henceforth, the Neumann-type boundary information of the source Fourier coefficients $v_n(\mathbf{x})$ for $1\le n \le N$ is given by
\begin{align}\label{g11}
	g_{1n}\left(\mathbf{x}\right)=\int_{\underline{k}}^{\overline{k}}k^{-2}\left[\frac{F_{1}\left(\mathbf{x},k\right)}{F_{0}\left(\mathbf{x},k\right)}-\frac{\partial_{z}u_{i}\left(\mathbf{x},k\right)}{u_{i}\left(\mathbf{x},k\right)}\right]\Psi_{n}\left(k\right)dk\quad\text{for }\mathbf{x}\in\Gamma.
\end{align}

Associating \eqref{mainsys1} with \eqref{g01}, \eqref{g02}, \eqref{g11} forms our system of coupled elliptic
equations, whose solution is the vector function $V\left( \mathbf{x}\right) $ that
contains all all of the Fourier coefficients $v_n$ for $1\le n\le N$.

\subsection{Convexified costs functional and theorems revisited}

It is evident that
\eqref{mainsys1} is a system of coupled quasi-linear elliptic equations. The nonlinear
terms are generated by products of gradients $\nabla v_{n}\left( \mathbf{x}\right) \cdot \nabla v_{m}\left( \mathbf{x}\right) $. Therefore, conventional least-squares methods, which minimize the differential functional, may not yield desirable results. Nonlinear problems often exhibit non-convex cost functionals, resulting in multiple local minima and ravines. Hence, a good initial guess must be chosen to reach the global minimizer.

To tackle nonlinear inverse problems, convexification is one of some numerical methods available. This method and its variants construct a suitable weighted cost functional that is strongly convex over a bounded set of a Hilbert space. With this approach, the existence and uniqueness of a minimizer can be proven without any restriction on the size of the set. Additionally, convergence towards the correct solution is guaranteed.

Introduce $\mu_{\lambda}=\mu_{\lambda}(z)=e^{-\lambda(R+r)^2}e^{\lambda(z-r)^2}$ as the Carleman Weight Function (CWF). Then, we consider the following weighted cost functional $J:[H^p(\Omega)]^{N}\to \mathbb{R}_{+}$, for $p> 3$,
\begin{equation}\label{JJ}
	J(V) = \sum_{l = 1}^N \int_{\Omega} \mu_\lambda^2\left|\sum_{n = 1}^NS_{ln}
	\Delta v_n + \sum_{n, m = 1}^N P_{lnm} \nabla v_i \cdot \nabla v_m + \sum_{n
		= 1}^N Q_{ln}\nabla v_n\right|^2 d\mathbf{x} + \varepsilon
	\|V\|^2_{[H^p(\Omega)]^N}.	
\end{equation}
Here, the CWF plays several important roles in the convexification of interest. First, the function helps to control the highly nonlinear terms in the target quasi-linear system. Second, the CWF appears to maximize the influence of the measured boundary data at $\Gamma$. Lastly, by the presence of such a function, one can prove that the cost functional is globally strongly convex.

From now onward, we state the minimization problem.\\
\textbf{Minimization problem.} Minimize the cost functional $J(V)$ on the set $\overline{B(M)}$,
\[
B\left(M\right)=\left\{ V\in\left[H^{p}\left(\Omega\right)\right]^{N}:\left\Vert V\right\Vert _{\left[H^{p}\left(\Omega\right)\right]^{N}}\le M\right\}.
\]

Now, we formulate theorems of our convergence results. The theorems were proven in, e.g., \cite{KhoaKlibanovLoc:SIAMImaging2020,LeNguyen:3DArxiv}. Therefore, their proofs are omitted.We begin
with the Carleman estimate for the continuous Laplacian.

\begin{theorem} 
	There exists constants $\lambda_{0}=\lambda_{0}\left(\Omega,r\right)\ge 1$ and $C=C\left(\Omega,r\right)>0$ such as for every $V\in H_{0}^{2}\left(\Omega\right)$ and for all $\lambda\ge\lambda_{0}$ the following Carleman estimate holds true:
	\begin{align*}
		\int_{\Omega}\mu_{\lambda}\left(z\right)\left|\Delta V\right|^{2}d\mathbf{x} & \ge\frac{C}{\lambda}\int_{\Omega}\mu_{\lambda}\left(z\right)\left(\left|\partial_{xx}V\right|^{2}+\left|\partial_{yy}V\right|^{2}+\left|\partial_{zz}V\right|^{2}+\left|\partial_{xy}V\right|^{2}+\left|\partial_{yz}V\right|^{2}+\left|\partial_{xz}V\right|^{2}\right)d\mathbf{x}\\
		& +C\lambda\int_{\Omega}\mu_{\lambda}\left(z\right)\left(\left|\nabla V\right|^{2}+\lambda^{2}\left|V\right|^{2}\right)d\mathbf{x}.
	\end{align*}
\end{theorem}

The next theorem is devoted to the global strong convexity of the cost functional $J(V)$.

\begin{theorem} \label{convex}
	The functional $J(V)$ defined in \eqref{JJ} has its Fr\'echet derivative $DJ$ for all $V\in \overline{B(M)}$. Moreover, we can find a sufficiently large $\lambda = \lambda(M,\Omega) >0$ such that  $J(V)$ is strongly convex on $\overline{B(M)}$. In particular, for all $V_2,V_1\in \overline{B(M)}$, we have
	\begin{align*}
		J\left(V_{2}\right)  -J\left(V_{1}\right)-DJ\left(V_{1}\right)\left(V_{2}-V_{1}\right)
	 \ge C\left\Vert V_{2}-V_{1}\right\Vert _{\left[H^{2}\left(\Omega\right)\right]^{N}}^{2}+\varepsilon\left\Vert V_{2}-V_{1}\right\Vert _{\left[H^{p}\left(\Omega\right)\right]^{N}}^{2},
	\end{align*}
	where $C=C(M,\Omega)>0$.
\end{theorem}

As a by-product of Theorem \ref{convex}, the existence and uniqueness of a minimizer $V_{\min}$ in $\overline{B(M)}$ are guaranteed. Moreover, we obtain the Lipschitz continuity of the Fr\'echet derivative $DJ$ on $\overline{B(M)}$; see e.g. \cite[Theorem 5.2]{KhoaKlibanovLoc:SIAMImaging2020} and some other references cited therein.

The convergence result follows from \cite[Theorem 6]{Poyak1987}. Following the Tikhonov regularization concept \cite{Beilina2012}, we assume the existence of the exact solution $V_{*}\in [H^{p}(\Omega)]^{N}$ of system  \eqref{mainsys1} and that it satisfies the noiseless data $g_{0}^{*}$ and $g_1^{*}$. Here, $g_{0}^{*}$ and $g_1^{*}$ are, respectively, corresponding to the noisy boundary data $g_0$ and $g_1$, whose elements are defined in \eqref{g01}, \eqref{g02}, \eqref{g11}.

Let $\delta > 0$ be the noise level. We assume that there exists an error function $\mathcal{E}\in [H^{p}(\Omega)]^{N}$ satisfying
\[
\begin{cases}
	\left\Vert \mathcal{E}\right\Vert _{\left[H^{p}\left(\Omega\right)\right]^{N}}\le\delta,\\
	g_{0}=g_{0}^{*}+\left.\partial_{\nu}\mathcal{E}\right|_{\partial\Omega},\\
	g_{1}=g_{1}^{*}+\left.\mathcal{E}\right|_{\Gamma}.
\end{cases}
\]
Next, we assume the existence of a function $\mathcal{V}$ such that $\partial_{\nu}\mathcal{V} = g_0$ on $\partial\Omega$ and $\mathcal{V} = g_1$ on $\Gamma$. Consider $V_{\min,\varepsilon,\delta}\left(\mathbf{x}\right)=V_{\min}\left(\mathbf{x}\right)+\mathcal{V}\left(\mathbf{x}\right)$ for $\mathbf{x}\in\Omega$. The convergence theorem is stated in the following.

\begin{theorem}\label{conv}
	Assume that
	\[
	\max\left\{ \left\Vert V_{*}\right\Vert _{\left[H^{p}\left(\Omega\right)\right]^{N}},\left\Vert \mathcal{V}\right\Vert _{\left[H^{p}\left(\Omega\right)\right]^{N}}\right\} \le\frac{M}{3}-\delta.
	\]
	Then we can find a constant $C=C\left(\Omega,r,M\right)$ such that
	the following estimate holds true
	\[
	\left\Vert V_{\min,\varepsilon,\delta}-V_{*}\right\Vert _{\left[H^{2}\left(\Omega\right)\right]^{N}}\le C\left(\sqrt{\varepsilon}\left\Vert V_{*}-\mathcal{V}\right\Vert _{\left[H^{p}\left(\Omega\right)\right]^{N}}+\delta\right).
	\]
\end{theorem}

Since smallness conditions are not
imposed on $M$, then the above convergence estimate confirms the \textit{global convergence} of the minimizer of the cost functional $J(V)$ to the exact solution.

It now remains to discuss how to find $V_{\min,\varepsilon,\delta}$ by the so-called gradient descent method. Let $\eta \in (0,1)$. The gradient descent method is given as follows:
\begin{align}\label{gd}
	V^{(n)} = V^{(n-1)} - \eta DJ(V^{(n-1)}),\quad n=1,2,\ldots,
\end{align}
where $V^{(n)}$ denotes the $n$th iteration for the approximation of the minimizer $V_{\min,\varepsilon,\delta}$. In \eqref{gd}, we use the starting point $V^{(0)}\in B(M)$ being an arbitrary point in that particular set. Recall that by Theorem \ref{convex}, we obtain that $V_{\min,\varepsilon,\delta}$ in $\overline{B(M)}$. Cf. \cite[Theorem 2]{LeNguyen:2020}, if we assume further that the ball centered at $V_{\min,\varepsilon,\delta}$ with the radius $\left\Vert V^{\left(0\right)}-V_{\min,\varepsilon,\delta}\right\Vert_{\left[H^{2}\left(\Omega\right)\right]^{N}}$ is contained in $B(M)$, then the distance between the $n$th iteration $V^{(n)}$ and the minimizer $V_{\min,\varepsilon,\delta}$ is controlled well by that radius $\left\Vert V^{\left(0\right)}-V_{\min,\varepsilon,\delta}\right\Vert_{\left[H^{2}\left(\Omega\right)\right]^{N}}$. In particular, we formulate the following theorem, while its proof is omitted.

\begin{theorem}\label{theoremlast}
	Let $V^{\left(0\right)}\in B\left(M\right)$ and $V_{\min,\varepsilon,\delta}\in B\left(M\right)$ be such that the ball centered at $V_{\min,\varepsilon,\delta}$ with the radius $\left\Vert V^{\left(0\right)}-V_{\min,\varepsilon,\delta}\right\Vert_{\left[H^{2}\left(\Omega\right)\right]^{N}}$ is contained in $B(M)$. Then there exists a sufficiently
	small number $\eta_{0}\in\left(0,1\right)$ such that $V^{\left(n\right)}\subset B\left(M\right)$
	for all $n=1,2,\ldots$ and for all $\eta\in (0,\eta_{0})$. Moreover, there exists a number $\varsigma=\varsigma\left(\eta\right)\in\left(0,1\right)$
	such that
	\[
	\left\Vert V^{\left(n\right)}-V_{\min,\varepsilon,\delta}\right\Vert _{\left[H^{2}\left(\Omega\right)\right]^{N}}\le\varsigma^{n}\left\Vert V^{\left(0\right)}-V_{\min,\varepsilon,\delta}\right\Vert _{\left[H^{2}\left(\Omega\right)\right]^{N}}.
	\]
\end{theorem}

By Theorems \ref{conv} and \ref{theoremlast}, we obtain the strong convergence of the sequence $\left\{V^{n}\right\}_{n=0}^{\infty}$ toward the exact solution $V_{*}$. Particularly, by the triangle inequality, it holds true that
\[
\left\Vert V^{\left(n\right)}-V_{*}\right\Vert _{\left[H^{2}\left(\Omega\right)\right]^{N}}\le C\left(\sqrt{\varepsilon}\left\Vert V_{*}-\mathcal{V}\right\Vert _{\left[H^{p}\left(\Omega\right)\right]^{N}}+\delta\right)+\varsigma^{n}\left\Vert V^{\left(n\right)}-V_{*}\right\Vert _{\left[H^{2}\left(\Omega\right)\right]^{N}}.
\]

\section{Numerical experiments}
\label{sec num}

The numerical results performed in this section are all obtained with experimental data. Those are data collected at the microwave facility of the University of North Carolina at Charlotte (UNCC), USA. For the sake of simplicity, we refer the reader to \cite{VoKlibanovNguyen:IP2020,Khoaelal:IPSE2020} for details of the experimental setup we establish at the University of North Carolina at Charlotte. Thereby, we only mention below key elements of our experimental configuration. Even though those publications \cite{VoKlibanovNguyen:IP2020,Khoaelal:IPSE2020} focus wholly on the CIP with multiple point sources and a fixed frequency, the set of data collected at that time is variable in both source locations and frequencies for trial-and-error. For each source position, our raw data set consists of back-scattering data corresponding to 201 frequency values uniformly distributed between 1 GHz to 10 GHz. Therefore, we are capable of using those data to verify the numerical performance of the convexification method for the CIP in the context of multiple frequencies and a fixed point source.

\subsection{Experimental configuration and computational settings}\label{subsec:setting}

The experiment conducted at UNCC involves practical data of five (5) experimental objects buried in a sandbox. Those tested objects were prepared to mimic explosive-like devices often seen in the battlefield. Typically, we classify them as metallic and non-metallic objects:
\begin{enumerate}
	\item An aluminum tube that mimics the NO-MZ 2B, a Vietnamese anti-personnel fragmentation mine; cf. \cite{Banks98}.
	\item A glass bottle filled with clear water that is a good fit of the Glassmine 43 in terms of the material; cf. \cite{Office}. Reconstructing the shape of the bottle is challenging as this object comes with a cap.
	\item An U-shaped piece of a dry wood that can be an example of the Schu-mine 42, a wood-based anti-personnel blast mine. Compared to the glass bottle above, this piece of dry wood has a very complicated non-convex shape.
	\item A metallic letter `A' that is to augment the complication in shape of metallic experimental object; compared to the aluminum cylinder.
	
	\item A metallic letter `O' that serves the same purpose as the letter `A'. It helps to test the numerical performance of the convexification method with varying levels of complexity in the shape of the object.
	
\end{enumerate}

In \cite{VoKlibanovNguyen:IP2020}, the last two tests (i.e. those with the metallic letters) were blinded. In this sense, we only knew their backscattering data
and that the experimental objects were buried close to the sand surface. The experimental results obtained in that publication, however, turn the blinded tests to be unblinded. Therefore, in the present paper, our numerical experiments are demonstrated with all unblinded tests.

For every test, the experimental object is placed  inside of a rectangular box filled with moisture-free sand, which is then referred to as a sandbox. This man-made sandbox is framed by some wood materials, and its back and front sides are covered by a 5-cm layer of Styrofoam. The front side we mean here is closer to the standard antenna, compared to the back one. In our configuration, the antenna plays a role in sending incident waves toward the sandbox. Then, there is a rectangular measurement surface of dimensions $100\times 100$ ($\text{cm}^2$) behind the antenna to collect the backscattering data. Experimentally, this surface is discretized in an equidistant mesh of 2-cm mesh-width. Moreover, the horizontal and vertical sides of this surface define, respectively, the $x$- and $y$-axes of our coordinate system and thus, the $z$-axis is the orthogonal one to our
measurement plane. As to the burial depth of the experimental object, it is a few centimeters away from the front Styrofoam. This is relevant to real-world applications that landmines are at most 10 (cm) away from the ground surface in the battlefield; cf. \cite{Daniels2006}.

In the sequel, we consider dimensionless variables as $\mathbf{x}' = \mathbf{x}/(10 \text{ cm})$ and for simplicity, we use the same notations as in the theoretical part. In this regard, the dimensions in our computations are 10 times less than the
real ones in centimeters. For instance, our $100\times 100$ ($\text{cm}^2$) measurement plane is understood as a $10\times 10$ surface in the dimensionless regime. Now, we introduce our computational setup in this dimensionless setting. According to our experiment, the
distance between the measurement surface and the sandbox with the front Styrofoam is 11.05. We also find that the length in the $z$ direction of our sandbox without the Styrofoam is about 4.4. As the Styrofoam layer is bent by the intensity of dry sand, we deliberately reduce 10\% of this length. All of these result in the choice of our computational domain $\Omega = \left\{\mathbf{x}\in \mathbb{R}^3: \left|x\right|,\left|y\right|<5,\left|z\right|<2\right\}$. In other words, we take $R=5$, $b=2$, and the center of the sandbox is taken as the origin of our coordinate system. As to the source position, in our numerical results below we choose the one adjacent right to the origin of the line of sources in \cite{VoKlibanovNguyen:IP2020,Khoaelal:IPSE2020}. The location of this source is $(0.1,0,-9)$.

Our raw data are measured far away from the sandbox. Cf. \cite{Liem2018,VoKlibanovNguyen:IP2020,Khoaelal:IPSE2020}, we observed numerically that these data lack quality due to many physical difficulties met in measurement process (antenna location, unwanted furniture with different materials, distracting signals). It will be then not good if we apply them directly to the minimization procedure. Thus, we employ the so-called data propagation technique to ``move'' these far-field data closer to the sandbox, which results in an approximation of the near-field data. It is also worth mentioning that the application of this data propagation procedure is helpful in reducing the size of the computational domain in the $z$-direction. Thus, it gives a better estimation of images of the experimental objects in $x,y$ coordinates. In this work, the near-field plane is chosen as $\Gamma = \left\{\mathbf{x}\in \mathbb{R}^3: \left|x\right|,\left|y\right|<5,\left|z\right|=2\right\}$ -- the front side of the sandbox.

\begin{figure}[!ht]
	\begin{center}
		\subfigure[Before frequencies filtering]{\label{ob1}
		\includegraphics[width=0.4\textwidth]{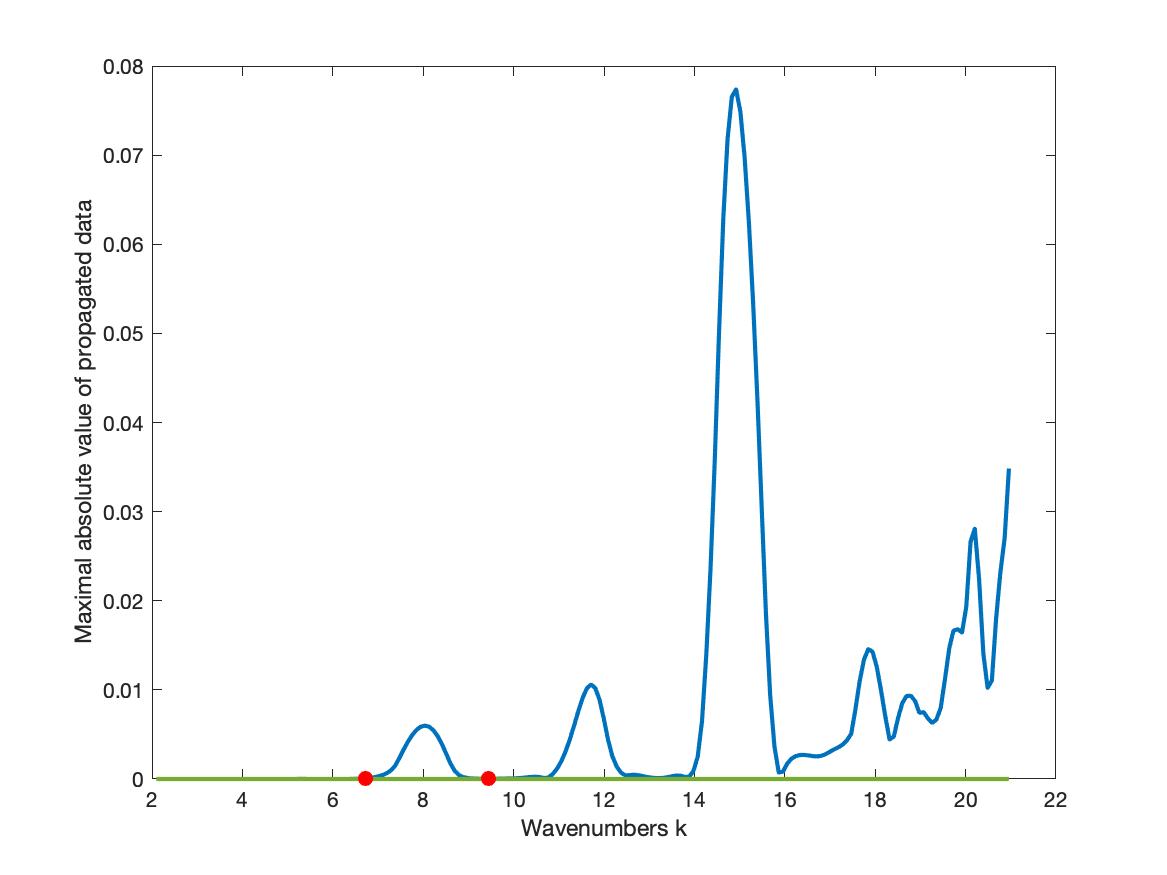} }
		\subfigure[After frequencies filtering]{\label{ob1_trunc}
		\includegraphics[width=0.4\textwidth]{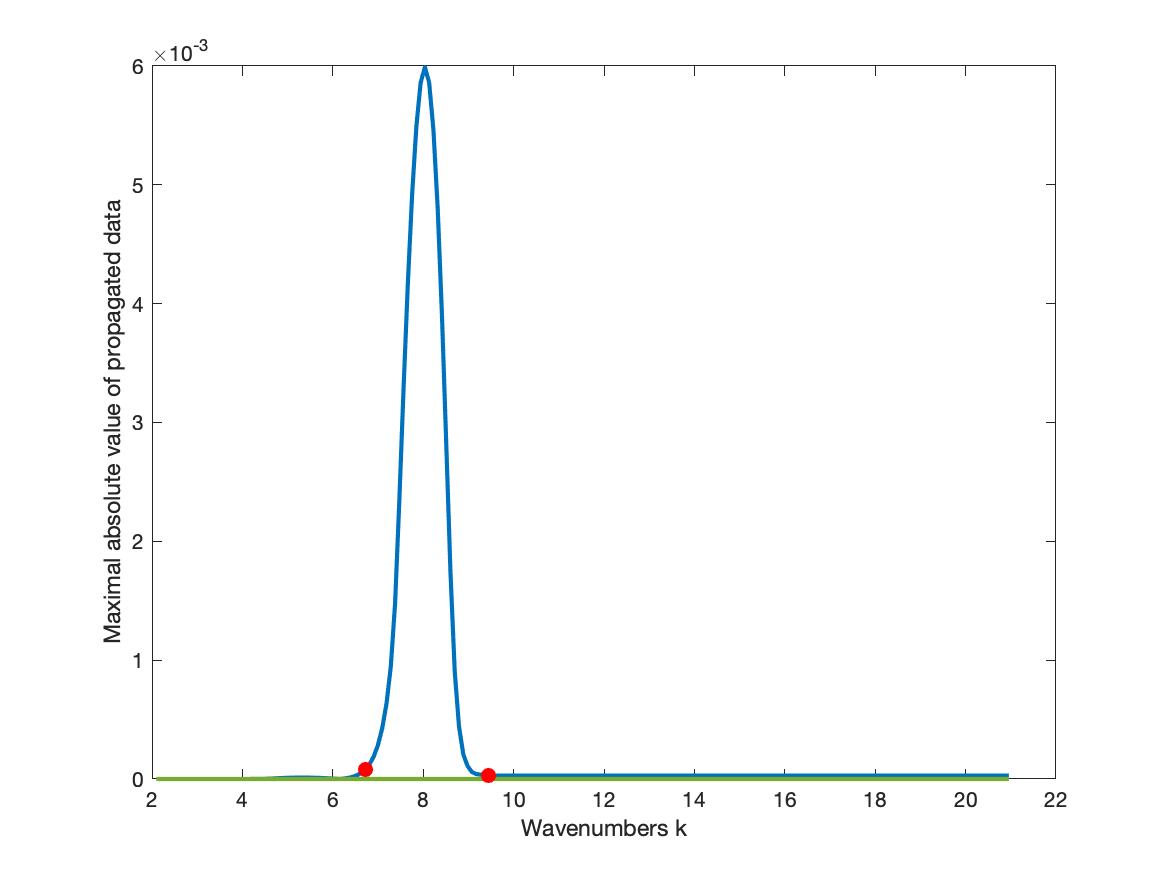}}
	\end{center}
	\caption{Illustration of the frequencies filtering for preprocessed data in Example 1. (a) The frequency dependent dynamics of the maximal absolute values of the experimental data after preprocessing procedure. All of these values are depicted for all wavenumbers $k \in [2.09,20.95]$ corresponding to the frequencies $\tilde{f} \in [1 \text{ GHz}, 10 \text{ GHz}]$. The red dots are imposed to indicate the wavenumber interval should be chosen. (b) The maximal absolute values of the processed data after frequencies filtering. \label{KU}}
	
\end{figure}

\begin{table}
	\begin{centering}
		\begin{tabular}{|c|c|c|c|c|c|}
			\hline 
			Example & 1 & 2 & 3 & 4 & 5\tabularnewline
			\hline 
			Object & Metallic cylinder & Bottle of water & Wood ``U'' & Metallic ``A'' & Metallic ``O''\tabularnewline
			\hline 
			Wavenumber $k$ & 6.72 - 9.45 & 5.87 - 8.60 & 18.22 - 20.96 & 10.68 - 13.41 & 8.70 - 11.43\tabularnewline
			\hline 
			Frequency (GHz) & 3.21 - 4.51 & 2.80 - 4.11 & 8.70 - 10 & 5.10 - 6.40 & 4.15 - 5.46\tabularnewline
			\hline 
		\end{tabular}
		\par\end{centering}
	\caption{The case of a
		single location of the source and multiple frequencies. Distinctive choices of the wavenumber
		interval for examples 1-5 and the corresponding frequencies. \label{tab:1}}
\end{table}

\begin{table}
	\begin{centering}
		\begin{tabular}{|c|c|c|c|c|c|}
			\hline 
			Example & 1 & 2 & 3 & 4 & 5\tabularnewline
			\hline 
			Object & Metallic cylinder & Bottle of water & Wood ``U'' & Metallic ``A'' & Metallic ``O''\tabularnewline
			\hline 
			Wavenumber $k$ & 8.51 & 6.62 & 11.43 & 9.55 & 8.79\tabularnewline
			\hline 
			Frequency (GHz) & 4.06 & 3.16 & 5.45 & 4.55 & 4.19\tabularnewline
			\hline 
		\end{tabular}
		\par\end{centering}
	\caption{The case of a
		single frequency and multiple sources. Distinctive choices of the wavenumber
		for examples 1-5 and the corresponding frequencies. \label{tab:2}}
\end{table}

We revisit the data propagation procedure that enables us to obtain the propagated data, termed as near-field data, from the raw data referred to as far-field data. We know that  $c(\mathbf{x}) = 1$ outside of the rectangular prism $\Omega = \left( -R,R\right) \times \left( -R,R\right) \times \left( -b,b\right)$ in $\mathbb{R}^{3}$. Therefore, the scattered wave $u_s$ in the half space $\left\{ z < -b\right\}$ satisfies the following system:
\begin{align}
	& \Delta u_s+k^{2} u_s = 0 \quad \text{in }\left\{ z < -b\right\},
	\label{eq:helms} \\
	& \lim_{r\rightarrow \infty }r\left( \partial _{r}u_s-\text{i}ku_s \right)
	=0\quad \text{for }r=\left\vert \mathbf{x}-\mathbf{x}_{\alpha }\right\vert ,%
	\text{i}=\sqrt{-1}.  \label{eq:somms}
\end{align}

As mentioned in section \ref{sec 2}, our experiments make use of the far-field data. Consequently, we have a dataset denoted as $u_s(x,y,-D)$, where $D > b$, while our objective of the data propagation procedure is to obtain $u_s(x,y,-b)$. Specifically, in our experiments, we have $D = 14$. The data propagation is obtained in the following way. Consider the Fourier transform of the scattered wave:
\[
\mathcal{F}\left(u_{s}\right)\left(\rho_{1},\rho_{2},z\right)=\frac{1}{2\pi}\int_{\mathbb{R}^{2}}u_{s}\left(\mathbf{x}\right)e^{-\text{i}\left(x\rho_{1}+y\rho_{2}\right)}dxdy,
\]
assuming that the corresponding integral is convergent. For $z<-b$, by applying the same Fourier transform to equation \eqref{eq:helms}, we find that
\begin{align*}
	\partial_{zz}^{2}\mathcal{F}\left(u_{s}\right)+\left(k^{2}-\rho_{1}^{2}-\rho_{2}^{2}\right)\mathcal{F}\left(u_{s}\right)=0.
\end{align*}
Solving the above differential equation gives the following relation between $\mathcal{F}\left(u_{s}\right)(z)$ at $z<-b$ and $\mathcal{F}\left(u_{s}\right)(z=-b)$, 
\begin{align}\label{dataprop1}
	\mathcal{F}\left(u_{s}\right)\left(z\right)=\begin{cases}
		\mathcal{F}\left(u_{s}\right)\left(-b\right)e^{\sqrt{\rho_{1}^{2}+\rho_{2}^{2}-k^{2}}\left(z+b\right)} & \text{for }\rho_{1}^{2}+\rho_{2}^{2}\ge k^{2},\\
		C_{1}e^{-\text{i}\sqrt{k^{2}-\rho_{1}^{2}-\rho_{2}^{2}}\left(z+b\right)}+C_{2}e^{\text{i}\sqrt{k^{2}-\rho_{1}^{2}-\rho_{2}^{2}}\left(z+b\right)} & \text{otherwise }.
	\end{cases}
\end{align}
Cf. \cite[Theorem 4.1]{Liem2018}, we can set $C_2 = 0$ in \eqref{dataprop1}. For $D$ relatively large, the value of the term in the first line of \eqref{dataprop1} is very small. Therefore, we can neglect the term with high frequencies and then use the inverse Fourier transform to obtain the near-field data:
\[
u_{s}\left(x,y,-b\right)=\frac{1}{\left(2\pi\right)^{2}}\int_{\rho_{1}^{2}+\rho_{2}^{2}<k^{2}}\mathcal{F}\left(u_{s}\right)\left(\rho_{1},\rho_{2},-D\right)e^{\text{i}\left[\sqrt{k^{2}-\rho_{1}^{2}-\rho_{2}^{2}}\left(-D+b\right)+x\rho_{1}+y\rho_{2}\right]}d\rho_{1}d\rho_{2}.
\]

As an slight improvement of the data propagation technique commenced in \cite{Liem2018}, we postulated a modified truncation procedure in \cite{VoKlibanovNguyen:IP2020} to remove possible random oscillations in the propagated data. The truncation procedure consists of two steps. First, we only preserve the propagated data whose values are at least 40 percents of the maximum absolute value. Then, we smooth those truncated data by the Gaussian filter. Observe that the smoothing process will average out the maximum value of the data, which may impact on the accuracy of the minimization result. In the second step, we add back some percents of the smoothed data to preserve the important ``peak'' that represents the maximum absolute value of the data. When doing so, we only need to multiply the smoothed data by a retrieval number computed by $\max\left(\left|\text{truncated data} \right|\right)/\max\left(\left|\text{smoothed data}\right|\right)$. This whole notion is mathematically specified in \cite{VoKlibanovNguyen:IP2020}, and the reader should be referred to that publication for any other details.

We now discuss the choice of an appropriate frequency interval since it does affect the quality of the frequency-dependent data applied to the minimization process. We remark that the raw data are frequency dependent in which the frequency unit is Hz (or $\text{s}^{-1}$). Cf. \cite{Liem2018,VoKlibanovNguyen:IP2020}, we formulate  the relation (in the dimensionless regime) between the wavenumber $k$ (with unit $\text{cm}^{-1}$) and the frequency, denoted by $\tilde{f}$, as follows:
\begin{align}\label{kless}
	k = \frac{2\pi \tilde{f}}{2997924580}.
\end{align} 
The choice of a frequency interval essentially relies on the performance of the data after preprocessing. Experimentally, it is different from one example to the others; see Table \ref{tab:1}. Following two criteria proposed in \cite{Liem2018}, we choose the wavenumber interval such that (1) the maximal absolute value of the processed data in this interval should not soar and plunge dramatically, and (2) for distinctive frequencies within this interval, these maxima at attained at the same coordinates (small deviation is acceptable) of the near-field plane. Once the interval is determined, we truncate all data that are outside of the chosen interval. Presenting the maximal absolute values of the experimental data in Example 1 after preprocessing for all wavenumbers $k$, Figure \ref{KU} exemplifies well the above-mentioned strategy. We find numerically that in Example 1, the interval of wavenumbers should be the vicinity of the first bump with a length of about 2.7 (see red dots in Figures \ref{ob1} and \ref{ob1_trunc}). Note that since our frequencies $\tilde{f}$ are between 1 GHz and 10 GHz, the corresponding wavenumber $k$ should range approximately from 2.09 to 20.95 using (\ref{kless}). We apply the same process to all other examples to choose appropriate wavenumber intervals for them. Tentatively, we call this process frequencies filtering.


\subsection{Minimization process}

Theoretically, our convexification method is globally convergent for any initial solution $V_0 \in \overline{B}$. However, to reduce the elapsed time of computations, we deliberately find the initial solution $V_0= \left(v_1^{(0)},\dots,v_N^{(0)}\right)^{\text{T}}$ that is close to $V$. Recall that $V=(v_1,v_2,\dots,v_N)^{\text{T}}$  is the solution of the following nonlinear elliptic system: 
\begin{equation}\label{mainsys}
	\sum_{n=1}^{N}S_{ln}\Delta v_{n}(\mathbf{x})+\sum_{n,m=1}^{N}P_{lnm}\nabla v_{n}(\mathbf{x})\cdot\nabla v_{m}(\mathbf{x})+\sum_{n=1}^{N}Q_{ln}(\mathbf{x})\cdot\nabla v_{n}(\mathbf{x})=0\quad\text{for all }l=\overline{1,N},\mathbf{x}\in\Omega,
\end{equation}
associated with the boundary conditions $\partial_{\nu}V|_\Gamma = g_0, V|_{\partial \Omega} = g_1$. Note that $g_0$ and $g_1$ are obtained from our experimental data after the frequencies filtering process. In (\ref{mainsys}), we indicate that for $i, j, l = \overline{1,N}$,
\begin{align*}
	& S_{ln}={\displaystyle \int_{\underline{k}}^{\overline{k}}\Psi_{i}^{\prime}(k)\Psi_{l}(k)dk,}\\
	& P_{lnm}={\displaystyle 2\int_{\underline{k}}^{\overline{k}}\Big(k^{2}\Psi_{n}(k)\Psi_{m}^{\prime}(k)+k\Psi_{n}(k)\Psi_{m}(k)\Big)\Psi_{l}(k)dk,}\\
	& Q_{ln}(\mathbf{x})={\displaystyle 2\int_{\underline{k}}^{\overline{k}}\left(\Psi_{n}^{\prime}(k)\frac{\nabla u_{i}(\mathbf{x},k)}{u_{i}(\mathbf{x},k)}+\Psi_{n}(k)\partial_{k}\frac{\nabla u_{i}(\mathbf{x},k)}{u_{i}(\mathbf{x},k)}\right)\Psi_{l}(k)dk},
\end{align*}
where the upper and lower bounds of $k$ are determined in the frequencies filtering mentioned above.

It is thus natural to take $V_0$ solutions to the corresponding linear system to (\ref{mainsys}). In this sense, we drop in (\ref{mainsys}) the nonlinear term containing $\nabla v_{n}\left( \mathbf{x}\right) \cdot $ $\nabla v_{m}\left( \mathbf{x}\right)$ and therefore, arrive at the following linear elliptic system:
\begin{equation}\label{initsys}
	\sum_{n=1}^{N}S_{ln}\Delta v_{n}(\mathbf{x})+\sum_{n=1}^{N}Q_{ln}(\mathbf{x})\cdot\nabla v_{n}(\mathbf{x})=0\quad\text{for all }l=\overline{1,N}
\end{equation}
associated with the same boundary conditions $g_0$ and $g_1$. By the natural linearity, system (\ref{initsys}) can be solved directly by the quasi-reversibility (QR) method involving the same Carleman weight function $\mu_{\lambda}=\mu_{\lambda}(z)=e^{-\lambda(R+r)^2}e^{\lambda(z-r)^2}$. In this regard, we minimize the following functional:
\begin{align}\label{QR}
\widetilde{J}\left(V\right)=\sum_{l=1}^{N}\int_{\Omega}\mu_{\lambda}^{2}\left|\sum_{n=1}^{N}S_{ln}\Delta v_{n}\left(\mathbf{x}\right)+\sum_{n=1}^{N}Q_{ln}\left(\mathbf{x}\right)\cdot\nabla v_{n}\left(\mathbf{x}\right)\right|^{2}d\mathbf{x} +\varepsilon\left\Vert V\right\Vert _{\left[H^{p}\left(\Omega\right)\right]^{N}}^{2}.
\end{align}
The solution $V_0$ obtained from solving (\ref{initsys}) will be used as the starting point of the minimization process.

Implementation of this QR method in a finite difference setting is detailed in \cite{LeNguyen:3DArxiv, LeNguyen:2020} and is analogous to the implementation of our cost functional $J(V)$ below.

As introduced in section \ref{sec:3}, the cost functional of our minimization process for system (\ref{mainsys}) is formulated as follows:
\begin{equation}
	J(V) = \sum_{l = 1}^N \int_{\Omega} \mu_\lambda^2\left|\sum_{n = 1}^NS_{ln}
	\Delta v_n + \sum_{n, m = 1}^N P_{lnm} \nabla v_i \cdot \nabla v_m + \sum_{n
		= 1}^N Q_{ln}\nabla v_n\right|^2 d\mathbf{x} + \varepsilon
	\|V\|^2_{[H^p(\Omega)]^N}.
	\label{J}
\end{equation}
Here, recall that  $\mu_{\lambda}=\mu_{\lambda}(z)=e^{-\lambda(R+r)^2}e^{\lambda(z-r)^2}$ involves the Carleman weight function $e^{\lambda(z-r)^2}$. In our numerical verification, we choose $\lambda=1.1$ and $r = 5.5$ in \eqref{Jdiscrete}. We remark that even though our theory is valid for sufficiently large values of $\lambda$, we have experienced numerically that we can choose a moderate value of $\lambda$ in [1,3]; see our previous works with both simulated and experimental data \cite{KhoaKlibanovLoc:SIAMImaging2020,Klibanov:2ndSAR2021,VoKlibanovNguyen:IP2020,Khoaelal:IPSE2020}. Below, we use the same value of $\varepsilon=10^{-9}$ for all examples. Also, instead of using a high regularity in the regularization term  $\varepsilon\|V\|_{[H^p(\Omega)]^N}$, we use only $\varepsilon\|V\|_{[H^2(\Omega)]^N}$. It reduces the computational complexity, while still providing a satisfactory numerical performance. As in \cite{VoKlibanovNguyen:IP2020}, the cut-off number  for our Fourier series is chosen as $N=6$ in all examples. Besides, the same parameters are used in our minmization of the quadratic functional \eqref{QR} of the QR method in all tests.

We now briefly mention the standard fully discrete version of the cost functional $J$ above. Let $N_x=N_y=51$ and $N_z = 21$ be the number of discrete points in $x,y$ and $z$ directions, respectively. Therefore, the same grid step size $h=0.2$ in these directions is used. For each $i=\overline{1,N}$, we denote by $v_i(x_{\mathfrak{i}}, y_{\mathfrak{j}}, z_{\mathfrak{l}})$ the corresponding discrete function of $v_i(x, y, z)$ defined at mesh-points $x_{\mathfrak{i}} = -R + \mathfrak{i}h, y_{\mathfrak{j}}=-R + \mathfrak{j}h, z_{\mathfrak{l}} = -b + \mathfrak{l}h$. Hereby, the corresponding Laplace operator in this finite difference setting is given by $\Delta^{h} = \partial_{xx}^{h} + \partial_{yy}^{h} + \partial_{zz}^{h}$, where, for interior grid points of $\Omega$, we consider
\begin{align*}
	\partial_{xx}^{h} v_{i}(x_{\mathfrak{i}}, y_{\mathfrak{j}}, z_{\mathfrak{l}}) = h^{-2} \left( v_{i}(x_{\mathfrak{i+1}}, y_{\mathfrak{j}}, z_{\mathfrak{l}}) - 2v_{i}(x_{\mathfrak{i}}, y_{\mathfrak{j}}, z_{\mathfrak{l}}) + v_{i}(x_{\mathfrak{i-1}}, y_{\mathfrak{j}}, z_{\mathfrak{l}}) \right),
\end{align*}
and the same is applied to the difference operators $\partial_{yy}^{h}, \partial_{zz}^{h}$. For the gradient operator, we consider $\nabla^{h} = (\partial_{x}^{h},\partial_{y}^{h},\partial_{z}^{h})$ with
\begin{align*}
	\partial_{x}^{h}v_{i}(x_{\mathfrak{i}}, y_{\mathfrak{j}}, z_{\mathfrak{l}}) = (2h)^{-1}(v_{i}(x_{\mathfrak{i+1}}, y_{\mathfrak{j}}, z_{\mathfrak{l}}) - v_{i}(x_{\mathfrak{i-1}}, y_{\mathfrak{j}}, z_{\mathfrak{l}})).
\end{align*}
Henceforth, the discrete version of $J$ corresponding to (\ref{J}) is given by
\begin{multline}
	J(V^{h}) = h^3\sum_{\mathfrak{i}=1}^{N_x}\sum_{\mathfrak{j}=1}^{N_y}\sum_{\mathfrak{l}=1}^{N_z}\sum_{l = 1}^N
	\mu_\lambda^2( z_{\mathfrak{l}})\Big|%
	\sum_{i = 1}^NS_{li} \Delta^{h} v_i(x_{\mathfrak{i}}, y_{\mathfrak{j}}, z_{%
		\mathfrak{l}}) \\
	+ \sum_{i, j = 1}^N P_{lij} \nabla^{h} v_i(x_{\mathfrak{i}}, y_{\mathfrak{j}},
	z_{\mathfrak{l}}) \cdot \nabla^{h} v_j(x_{\mathfrak{i}}, y_{\mathfrak{j}}, z_{%
		\mathfrak{l}}) + \sum_{i = 1}^N Q_{li}(x_{\mathfrak{i}}, y_{\mathfrak{j}},
	z_{\mathfrak{l}})\nabla^{h} v_i(x_{\mathfrak{i}}, y_{\mathfrak{j}}, z_{\mathfrak{%
			l}}) \Big|^2 \\
	+ \varepsilon h^3 \sum_{\mathfrak{i}=1}^{N_x}\sum_{\mathfrak{j}=1}^{N_y}\sum_{\mathfrak{l}=1}^{N_z} \sum_{l = 1}^N %
	\Big(|v_l(x_{\mathfrak{i}}, y_{\mathfrak{j}}, z_{\mathfrak{l}})|^2 + |\nabla^h
	v_l(x_{\mathfrak{i}}, y_{\mathfrak{j}}, z_{\mathfrak{l}})|^2 + |\Delta^{h}
	v_l(x_{\mathfrak{i}}, y_{\mathfrak{j}}, z_{\mathfrak{l}})|^2 \Big).
	\label{Jdiscrete}
\end{multline}

To speed up the computation process, we compute the gradient $DJ$ of the discrete functional $J$ in (\ref{Jdiscrete}) using the technique of Kronecker deltas; see in \cite{Kuzhuget:AA2010}. For brevity, we do not provide such a long formulation for the gradient $DJ$ here. Overall, the procedure to compute the approximate minimizer, denoted by $V^{\text{comp}}$, is described in Algorithm \ref{alg 1}. For the step size $\eta$ in Algorithm \ref{alg 1}, we briefly report that we start from $\eta = 10^{-1}$, and for each iterative step, if the value of the functional exceeds its value on the previous step, we replace the current step size $\eta$ by $\eta/2$. Otherwise, we keep it the same. We stop the minimization process via the gradient descent method when $\eta = 10^{-9}$.

\begin{algorithm}
	\caption{\label{alg 1} A numerical method to solve \eqref{mainsys}}
	\begin{algorithmic}[1]
		\State \label{choice} Choose a threshold error $\varepsilon > 0$.  
		\State  \label{choose U0} Set $m = 0$ and find an initial solution $V_0$ by solving \eqref{initsys}. 
		\State \label{step backward} Compute $V_{m + 1}$ using the gradient descent method for some step size $0 < \eta \ll 1$. 
		\State \label{step 4} If $\|V_{m + 1} - V_{m}\|_{[H^2(\Omega)]^{N}} < \varepsilon$, move forward to Step 5. Otherwise, set $m = m + 1$ and return Step 3.
		\State \label{step forward} Set $V^{\rm comp} = V_{m + 1}$
	\end{algorithmic}
\end{algorithm}

After obtaining the computed $N$-dimensional vector function $V^{\rm comp}$, we plug its components in the Fourier series that approximates $v$. Then, we compute the unknown dielectric constant by the following discrete formulation:
\begin{align*}
	& c^{h}(x_{\mathfrak{i}},y_{\mathfrak{j}},z_{%
		\mathfrak{l}})\\
	& =\text{mean}_{k}\left|-\text{Re}\left(\Delta^{h}v(x_{\mathfrak{i}},y_{\mathfrak{j}},z_{\mathfrak{l}},k)+k^{2}(\nabla^{h}v(x_{\mathfrak{i}},y_{\mathfrak{j}},z_{\mathfrak{l}},k))^{2}+\frac{2\nabla^{h}v(x_{\mathfrak{i}},y_{\mathfrak{j}},z_{\mathfrak{l}},k)\cdot\nabla^{h}u_{0}(x_{\mathfrak{i}},y_{\mathfrak{j}},z_{\mathfrak{l}},k)}{u_{0}(x_{\mathfrak{i}},y_{\mathfrak{j}},z_{\mathfrak{l}},k)}\right)\right|+1.
\end{align*}

To visually represent the reconstructed inclusion in each example, we use the isovalue function in MATLAB to generate 3D images. In cases where the inclusion possesses a high dielectric constant ($\ge 10$), we choose an isovalue of 20\%. Conversely, for inclusions with low dielectric constants ($<10$), we select an isovalue of 10\%.

\subsection{Numerical results}

\begin{figure}[!ht]
	\begin{center}
		\subfigure[Metallic cylinder]{\includegraphics[width=0.2\textwidth]{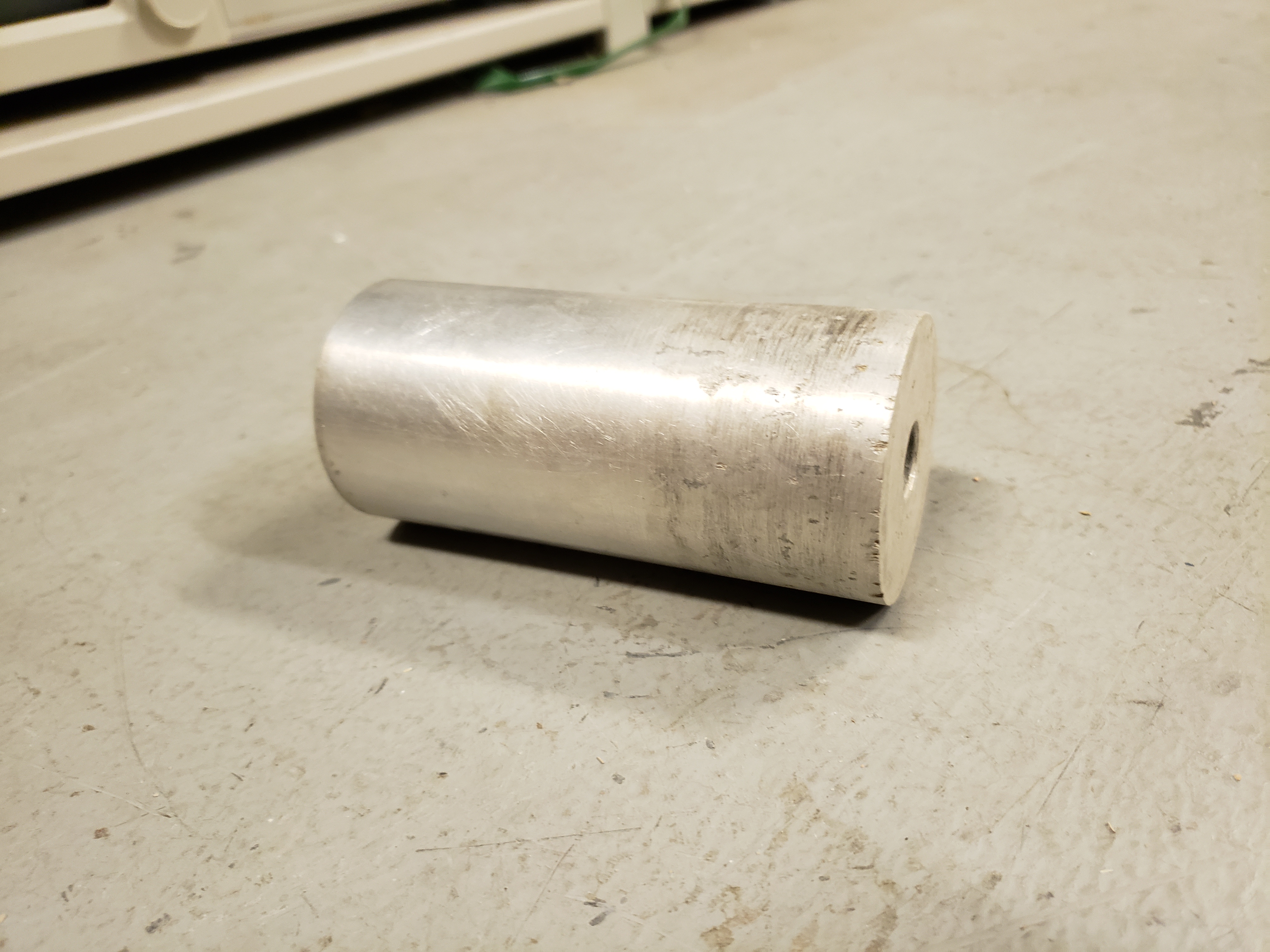}}

		\subfigure[The reconstructed metallic cylinder with one source and multiple frequencies]{\includegraphics[width=0.35\textwidth]{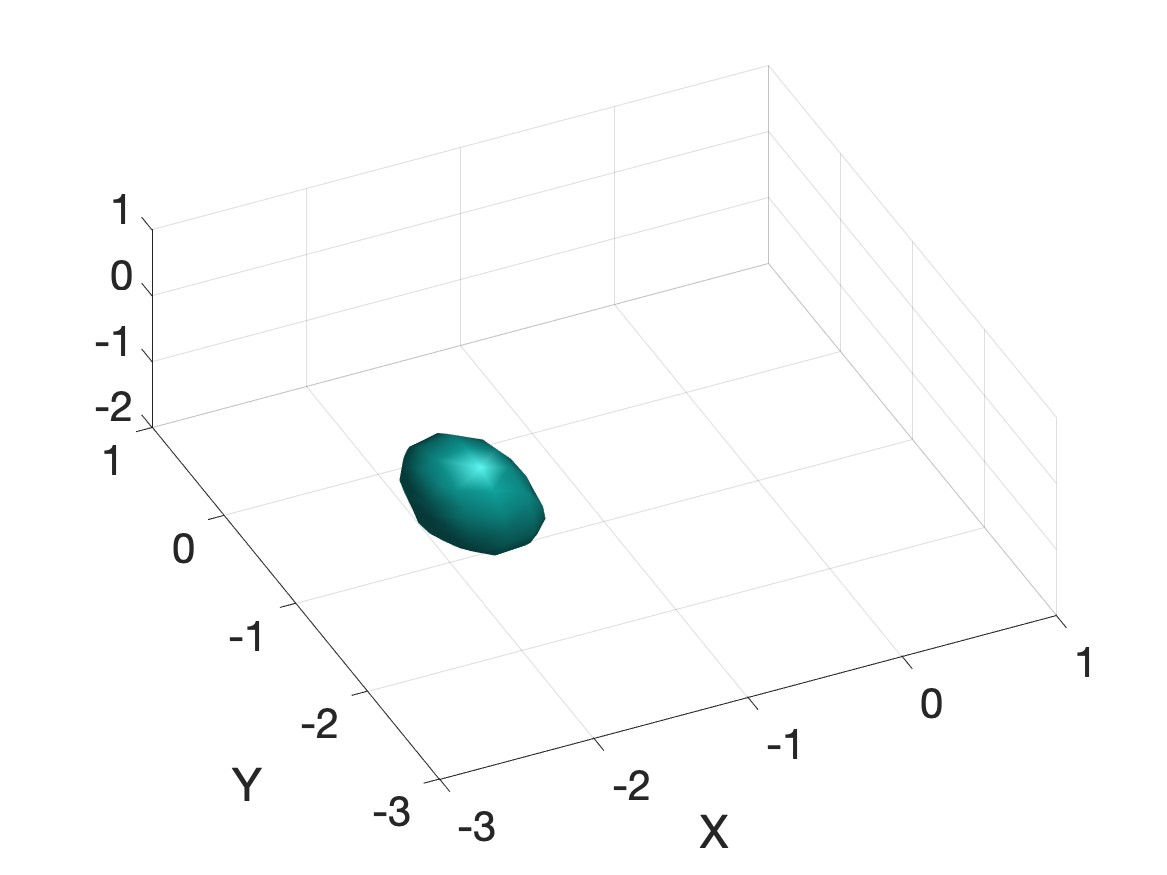}}

		\subfigure[The reconstructed of a metallic cylinder with multiple sources and a fixed frequency]{\includegraphics[width=0.35\textwidth]{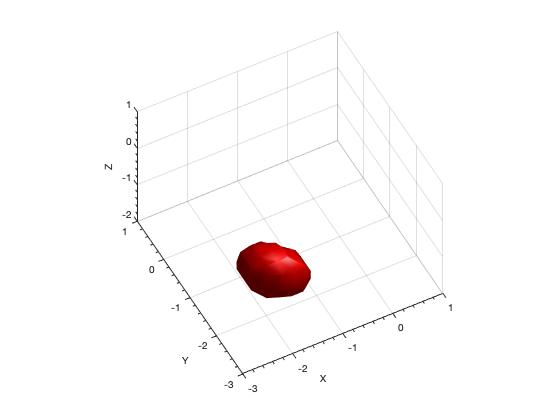}}
	\end{center}
	\caption{Example 1. Metallic cylinder. (a) Real photo of aluminum cylinder. (b) The reconstructed dielectric constant function by the first method with a fixed source and multiple frequencies. (c) The reconstructed dielectric constant function by the second method with multiple sources and a fixed frequency.\label{fig test 1}}
	
\end{figure}

\begin{figure}[!ht]
	\begin{center}
		\subfigure[Bottle of water]{
		\includegraphics[width=0.2\textwidth]{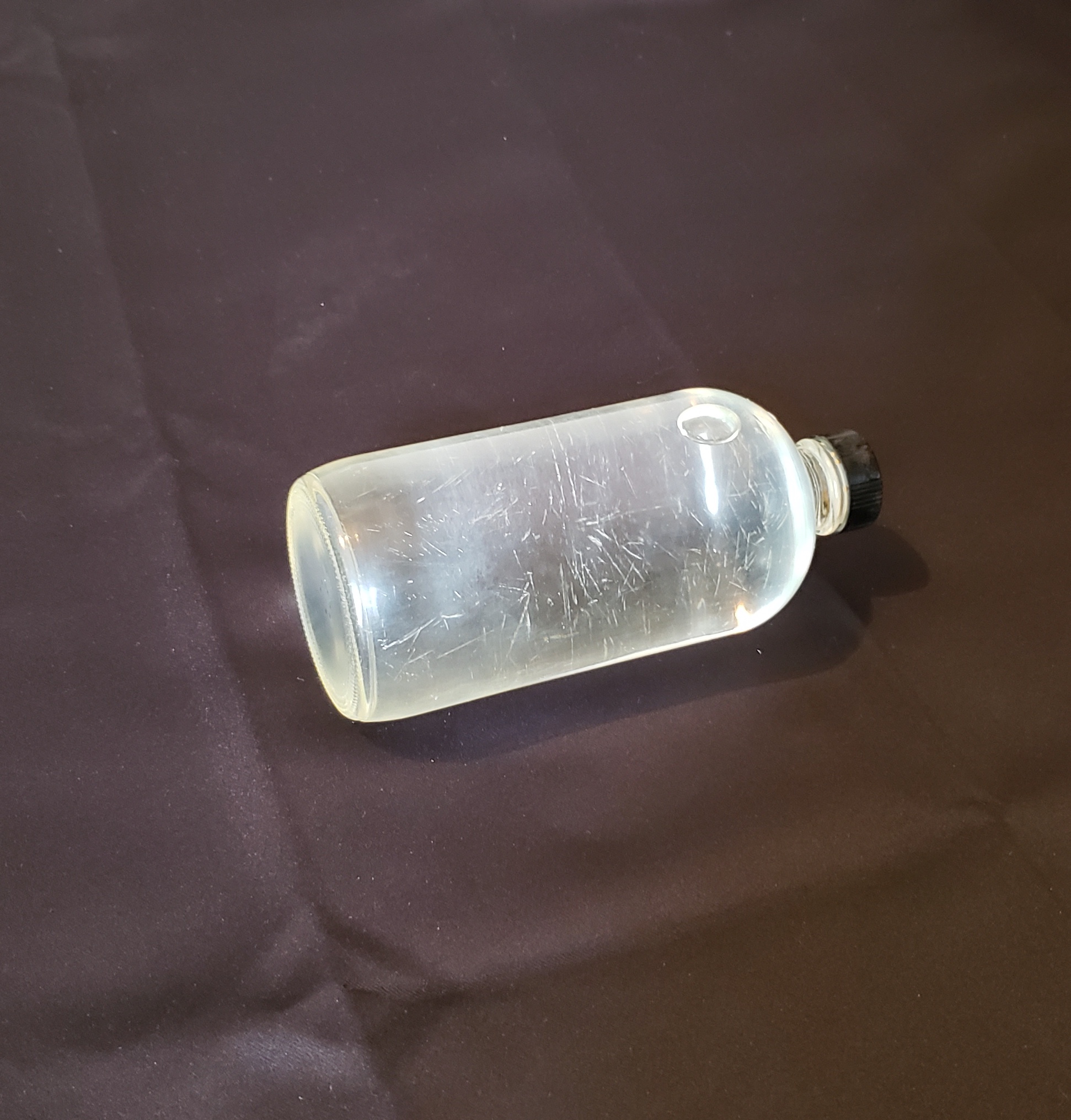}}
		
		\subfigure[The reconstructed bottle of water with one source and multiple frequencies]{\label{fig test 2b}
		\includegraphics[width=0.35\textwidth]{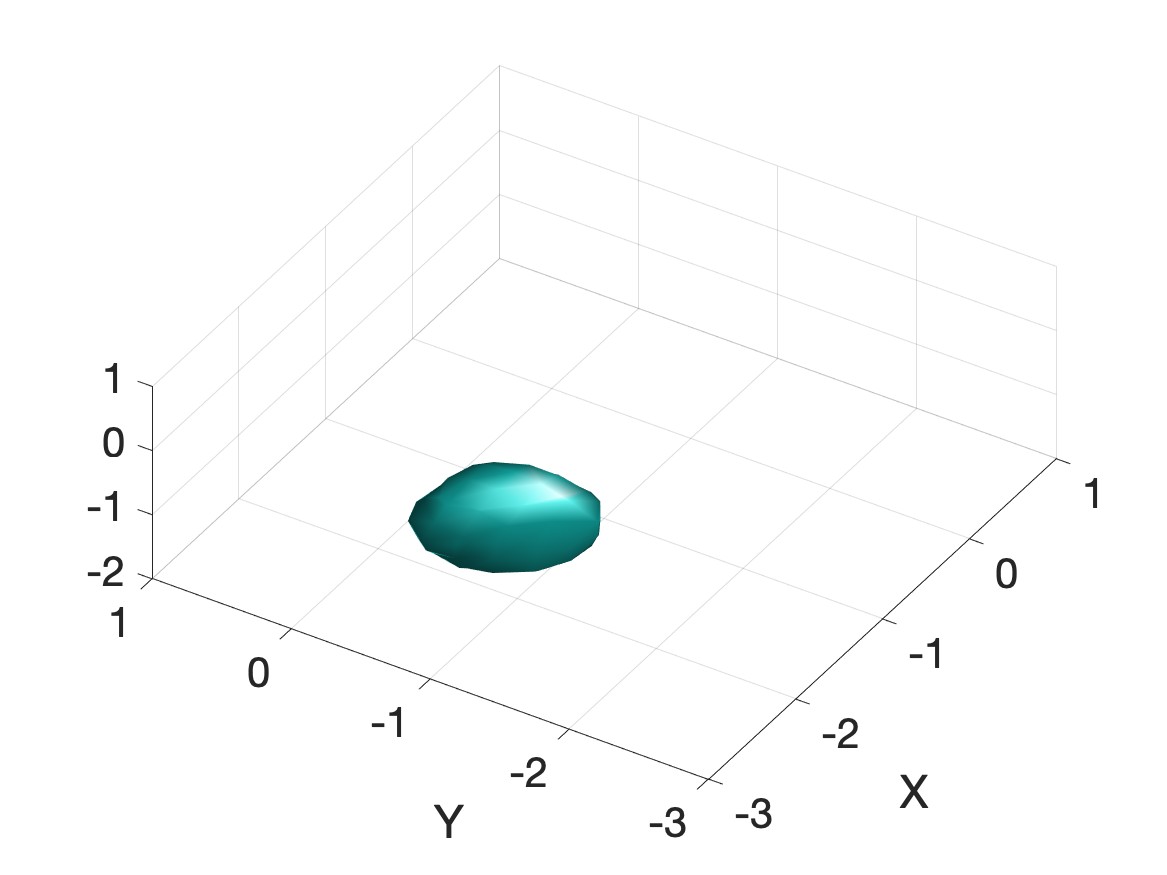}}
		
		\subfigure[The reconstructed bottle of water with multiple sources and a fixed frequency]{\label{fig test 2c}
		\includegraphics[width=0.35\textwidth]{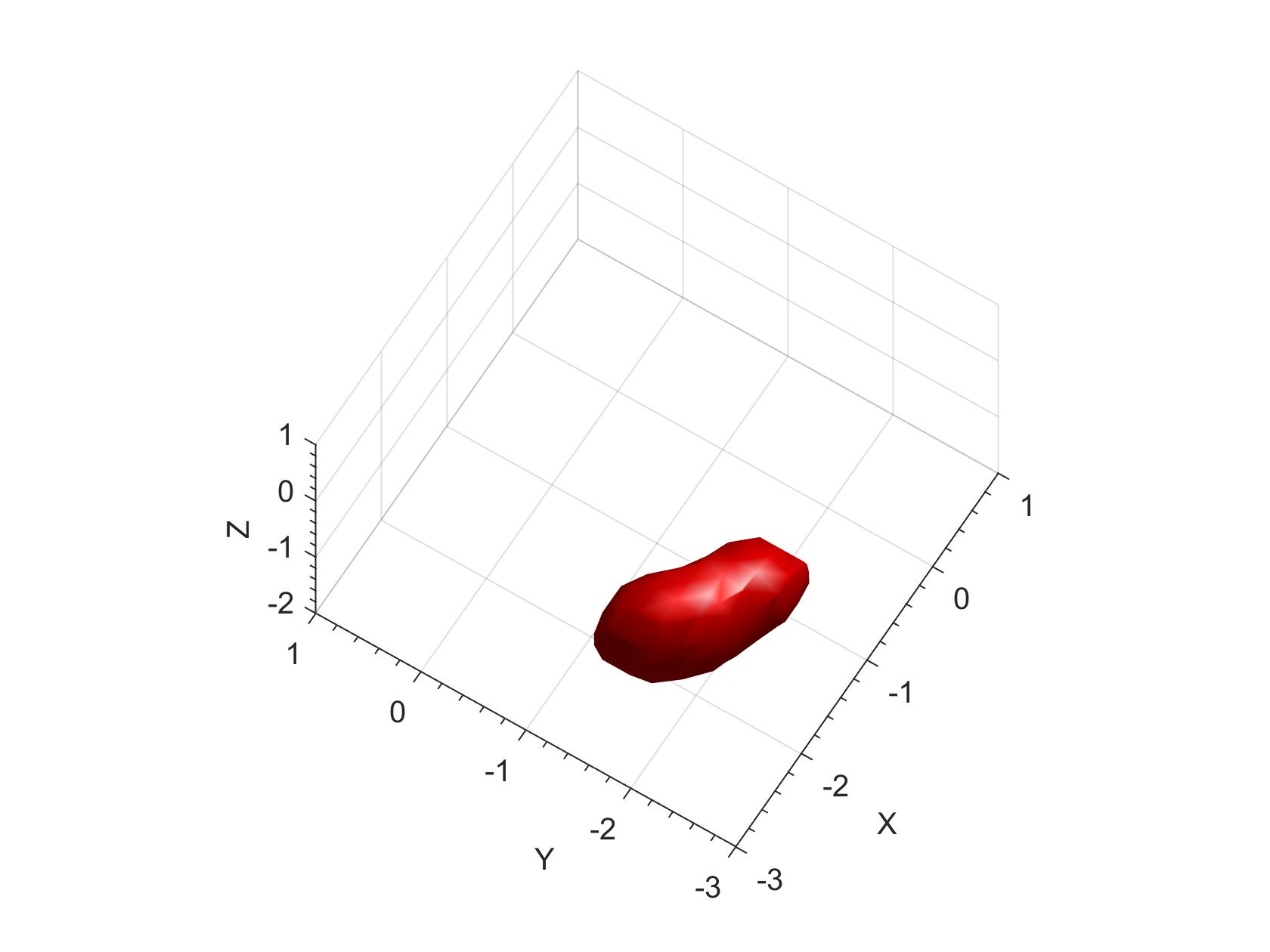}}
	\end{center}
	\caption{Example 2. The bottle of water. (a) The real image of the glass bottle of water. (b) The reconstructed dielectric constant function by the first method with a fixed source and multiple frequencies. (c) The reconstructed dielectric constant function by the second method with multiple sources and a fixed frequency. Clearly, we successfully detected the shape of the bottle of water by these two methods.	\label{fig test 2}}

\end{figure}

\begin{figure}[!ht]
	\begin{center}
		\subfigure[Wooden letter `U']{
		\includegraphics[width=0.2\textwidth]{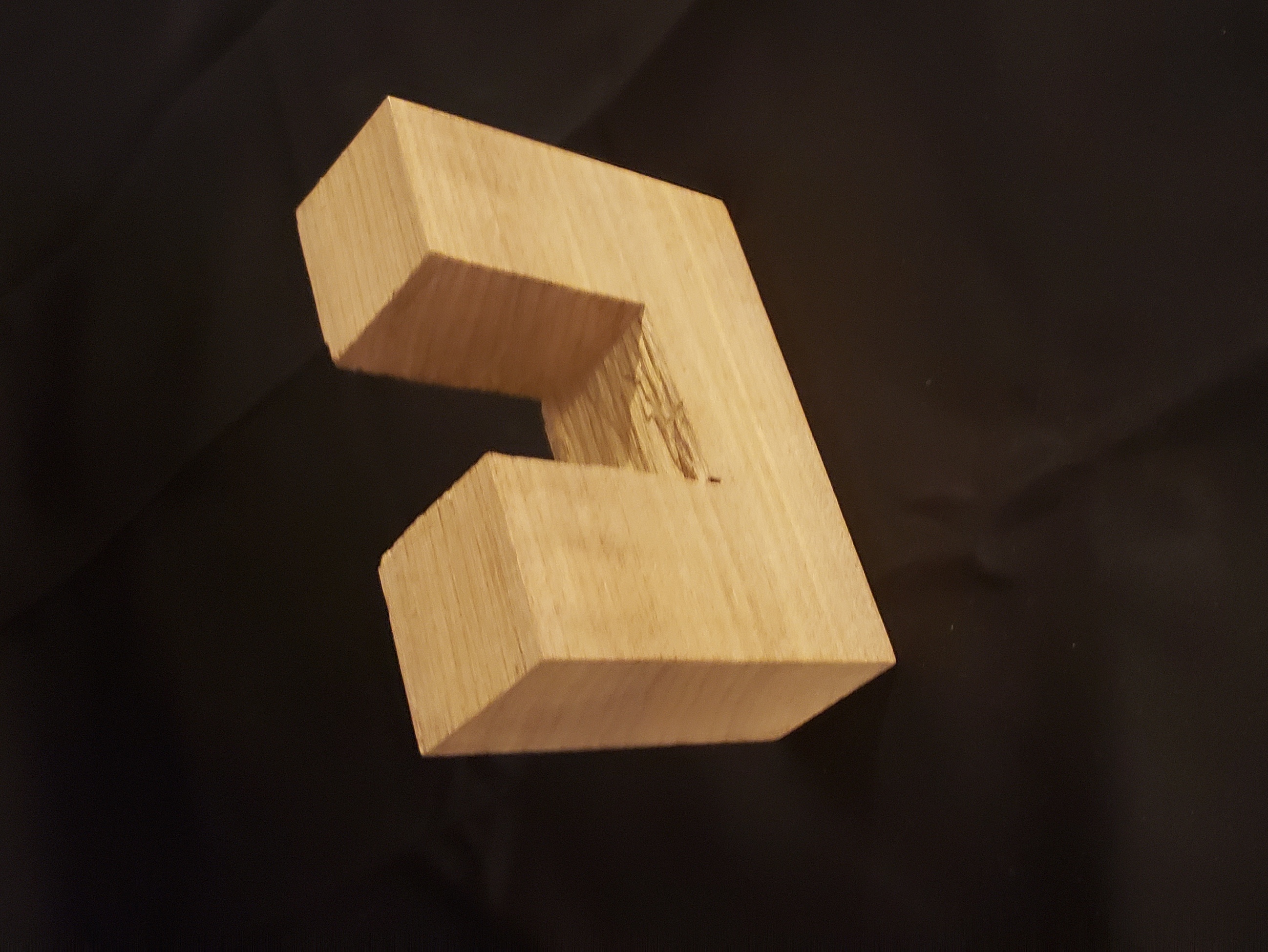}}
		
		\subfigure[The reconstructed letter `U' with one source and multiple frequencies]{
		\includegraphics[width=0.35\textwidth]{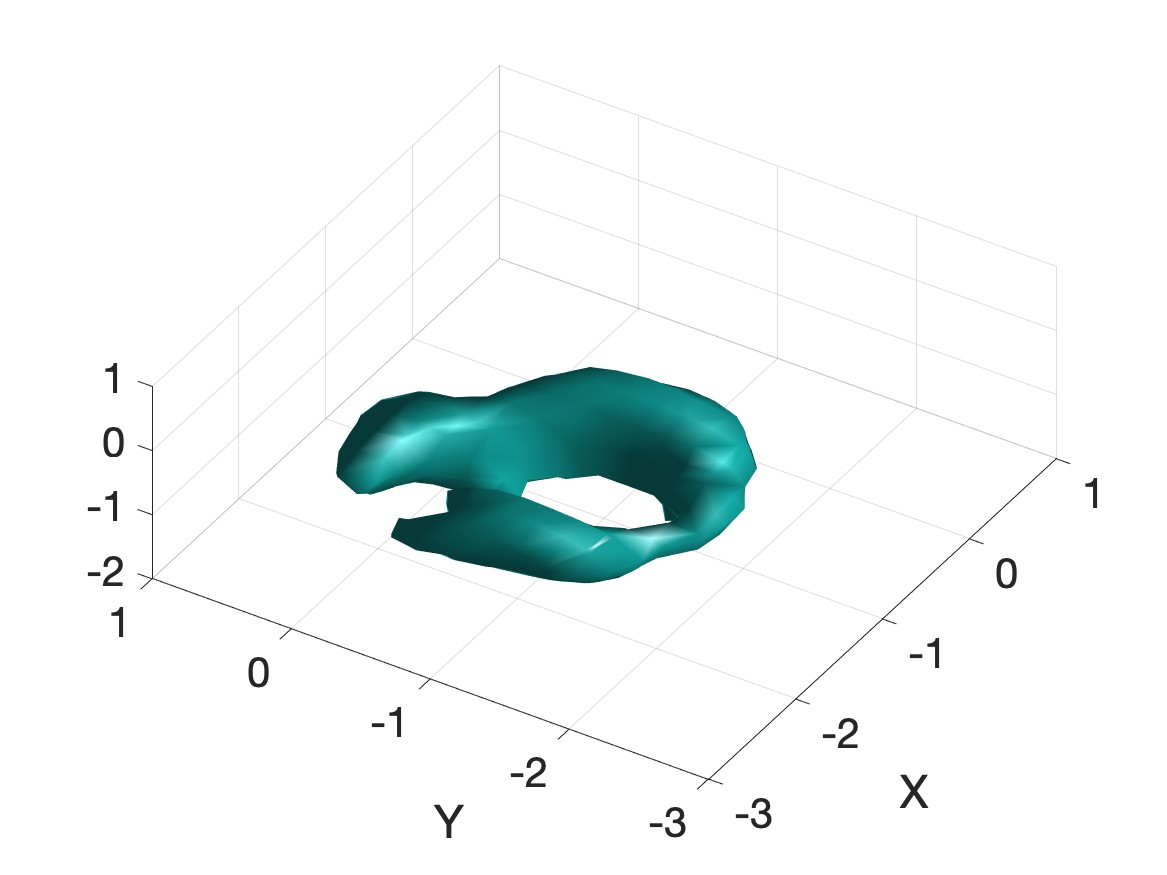}}
		
		\subfigure[The reconstructed letter `U' with multiple sources and a fixed frequency]{
		\includegraphics[width=0.35\textwidth]{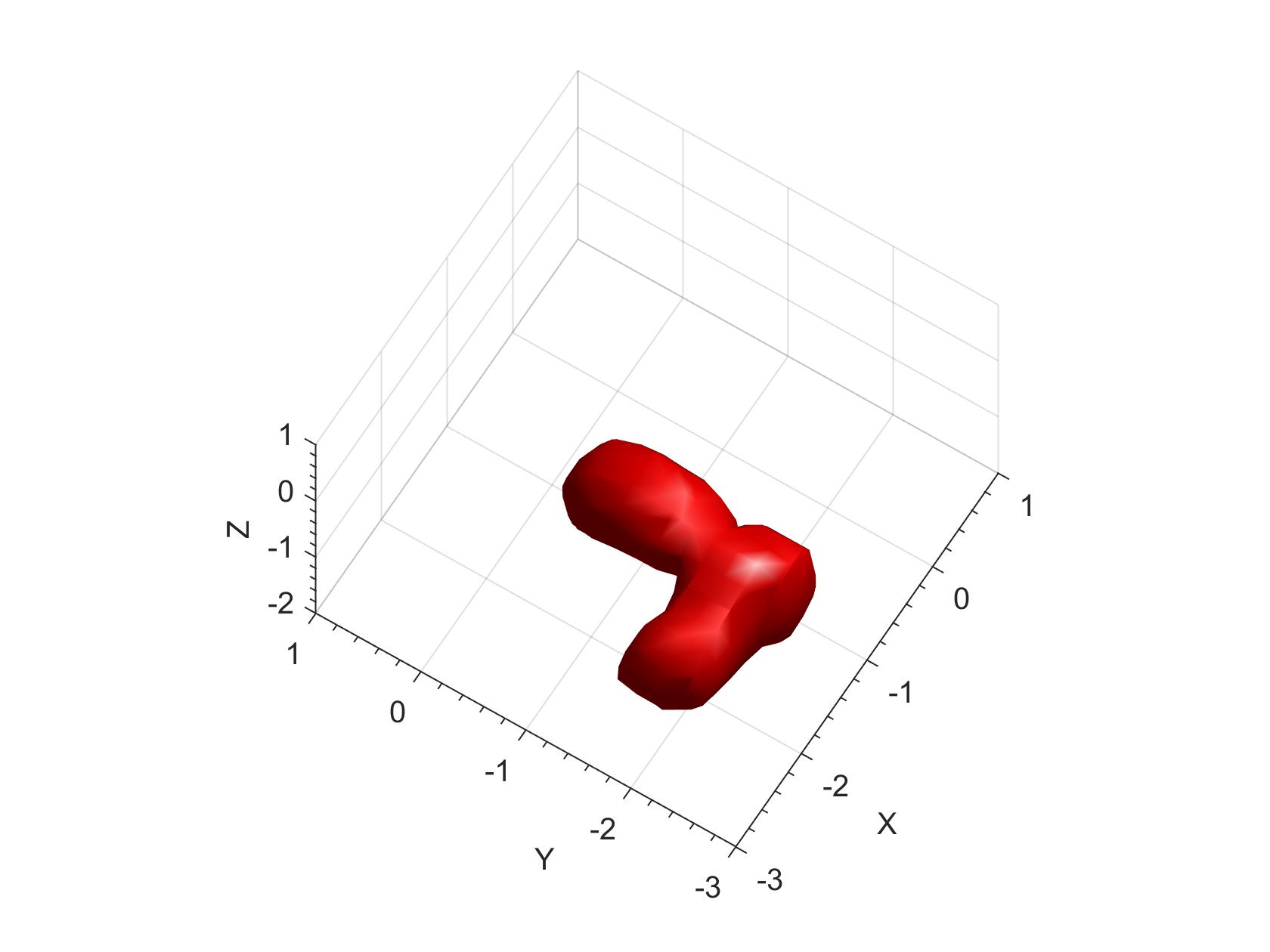}}
	\end{center}
	\caption{Example 3. The wooden letter `U'. (a) The real image of the U-shaped piece of dry wood. (b) The reconstructed dielectric constant function by the first method with a fixed source and multiple frequencies. (c) The reconstructed dielectric constant function by the second method with multiple sources and a fixed frequency. Note that the U shape can be seen clearly by the first method. It is well known that detecting non-convex objects with voids inside them is challenging but the first method can produce clearly letter U and the void inside it.\label{fig test 3}}
	
\end{figure}

\begin{figure}[!ht]
	\begin{center}
		\subfigure[Metallic letter `A']{
		\includegraphics[width=0.2\textwidth]{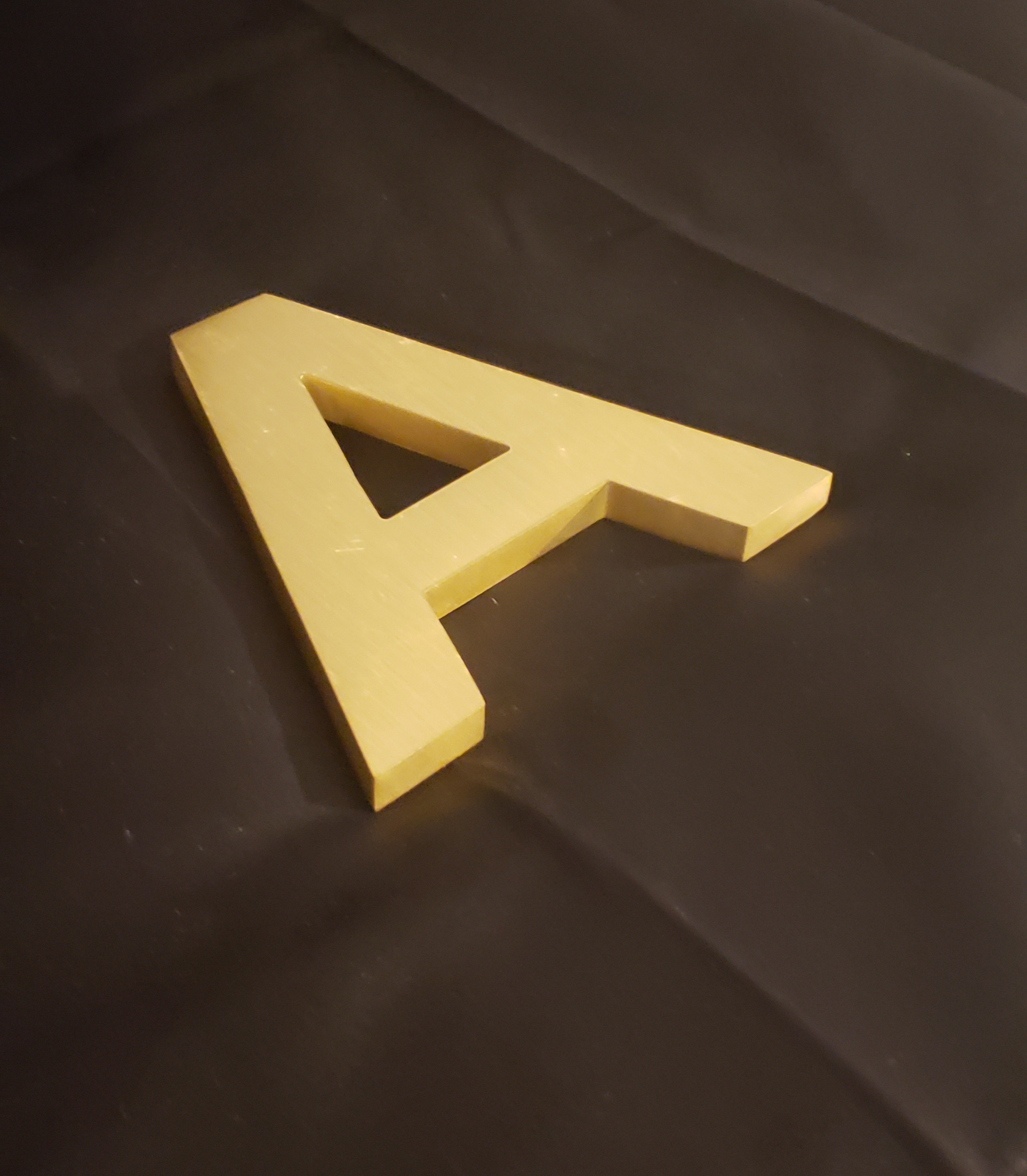}}
		
		\subfigure[The reconstructed letter `A' with one source and multiple frequencies]{
		\includegraphics[width=0.35\textwidth]{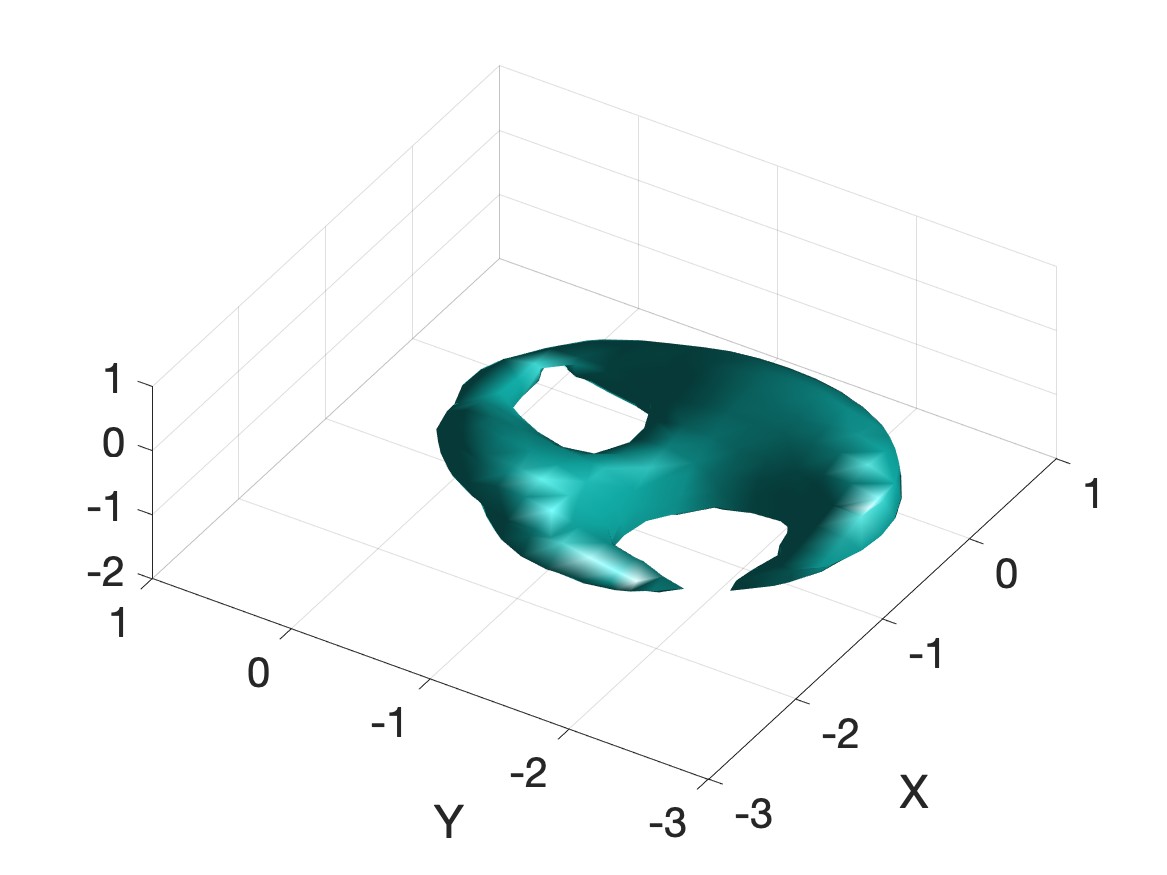}}
		
		\subfigure[The reconstructed letter `A' with multiple sources and a fixed frequency]{
		\includegraphics[width=0.35\textwidth]{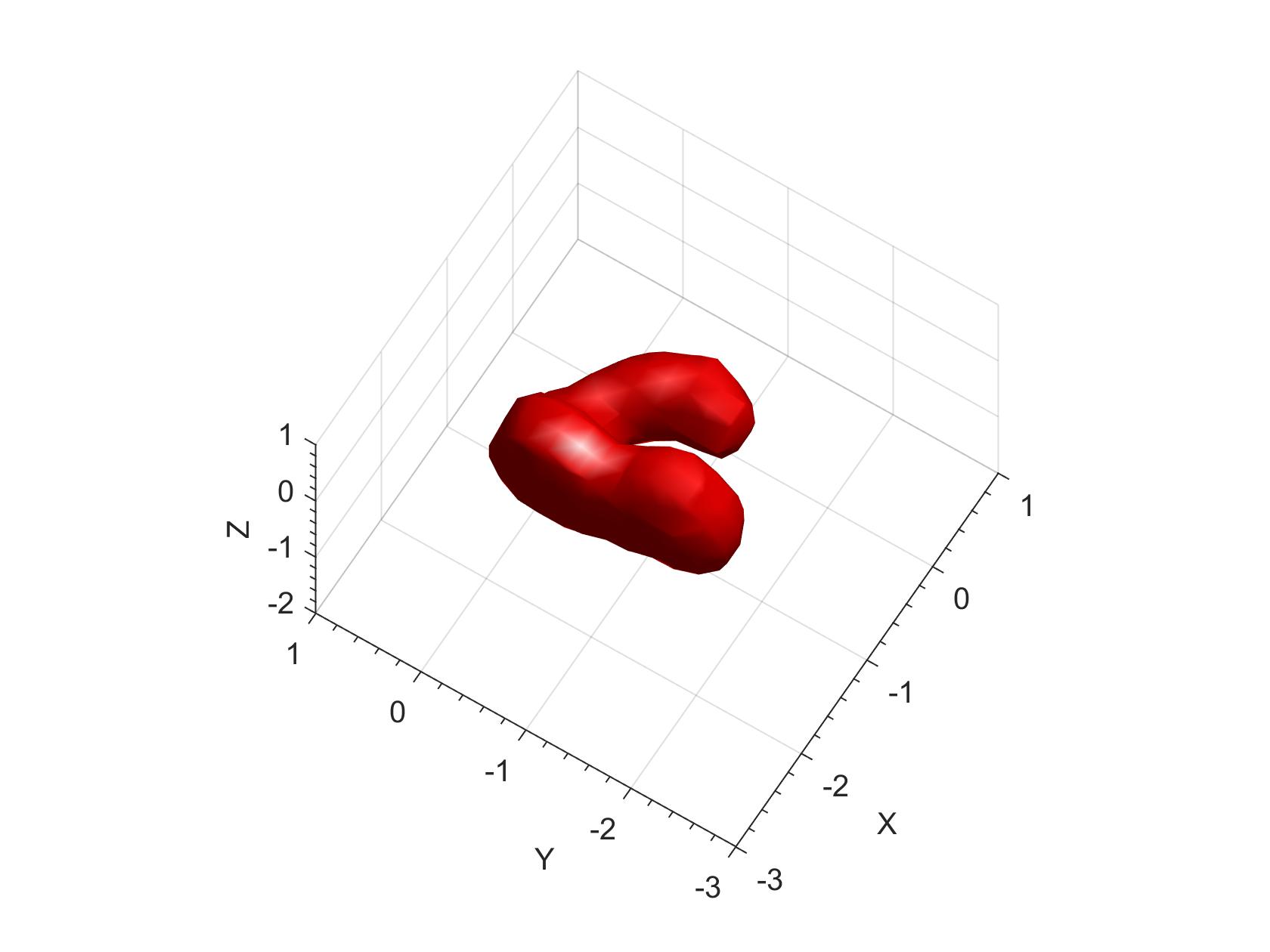}}
	\end{center}
	\caption{Example 4. The metallic letter `A'. (a) The real image of the A-shaped piece of metal. (b) The reconstructed dielectric constant function by the first method with a fixed source and multiple frequencies. (c) The reconstructed dielectric constant function by the second method with multiple sources and a fixed frequency. Note that letter A is produced perfectly with the void inside by the first method. It is much clearer than the result obtain by the second method.\label{fig test 4}}
	
\end{figure}

\begin{figure}[!ht]
	\begin{center}
		\subfigure[Metallic letter `O']{
		\includegraphics[width=0.2\textwidth]{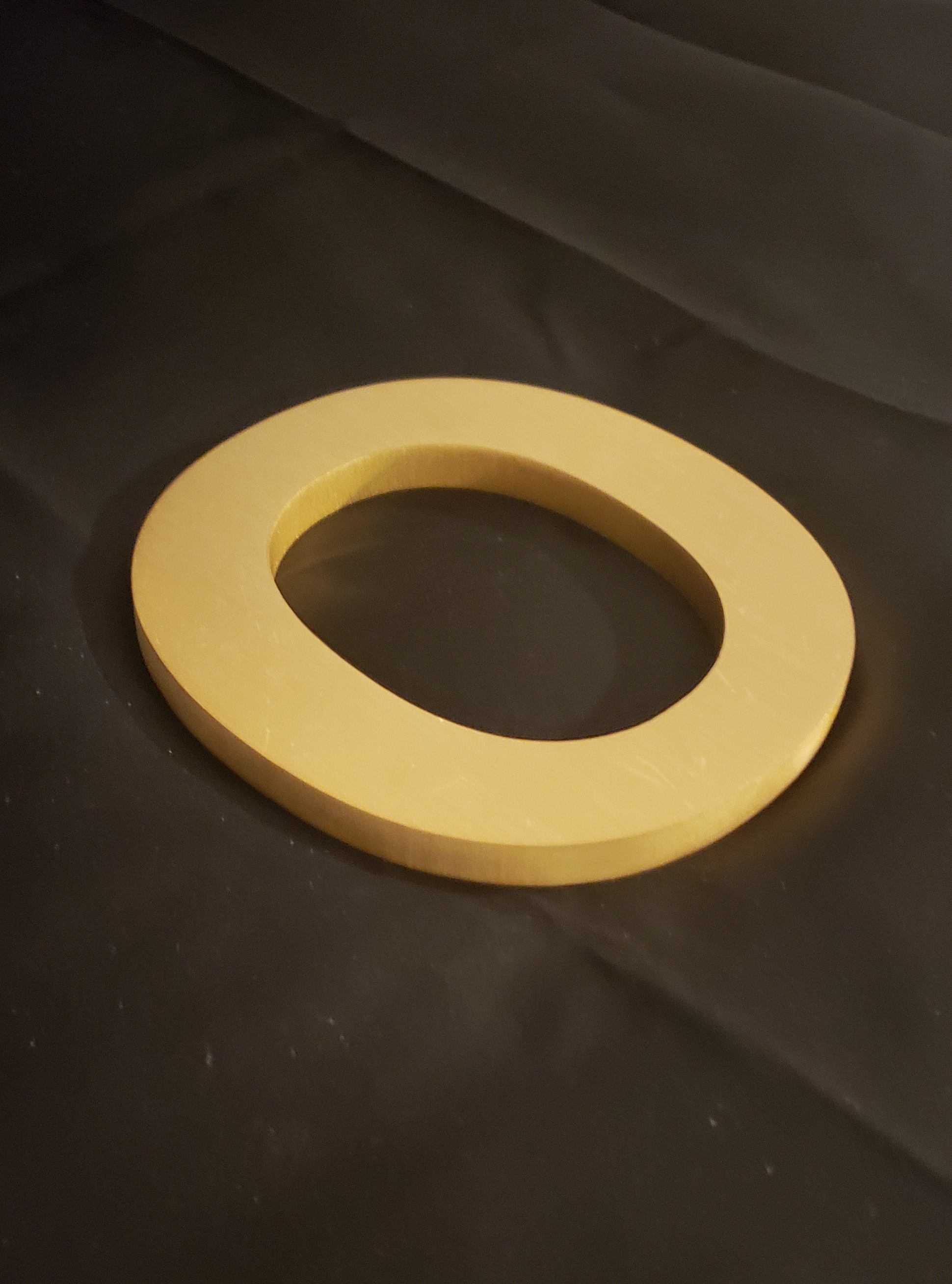}}
		
		\subfigure[The reconstructed letter
		`O' with one source and multiple frequencies]{
		\includegraphics[width=0.35\textwidth]{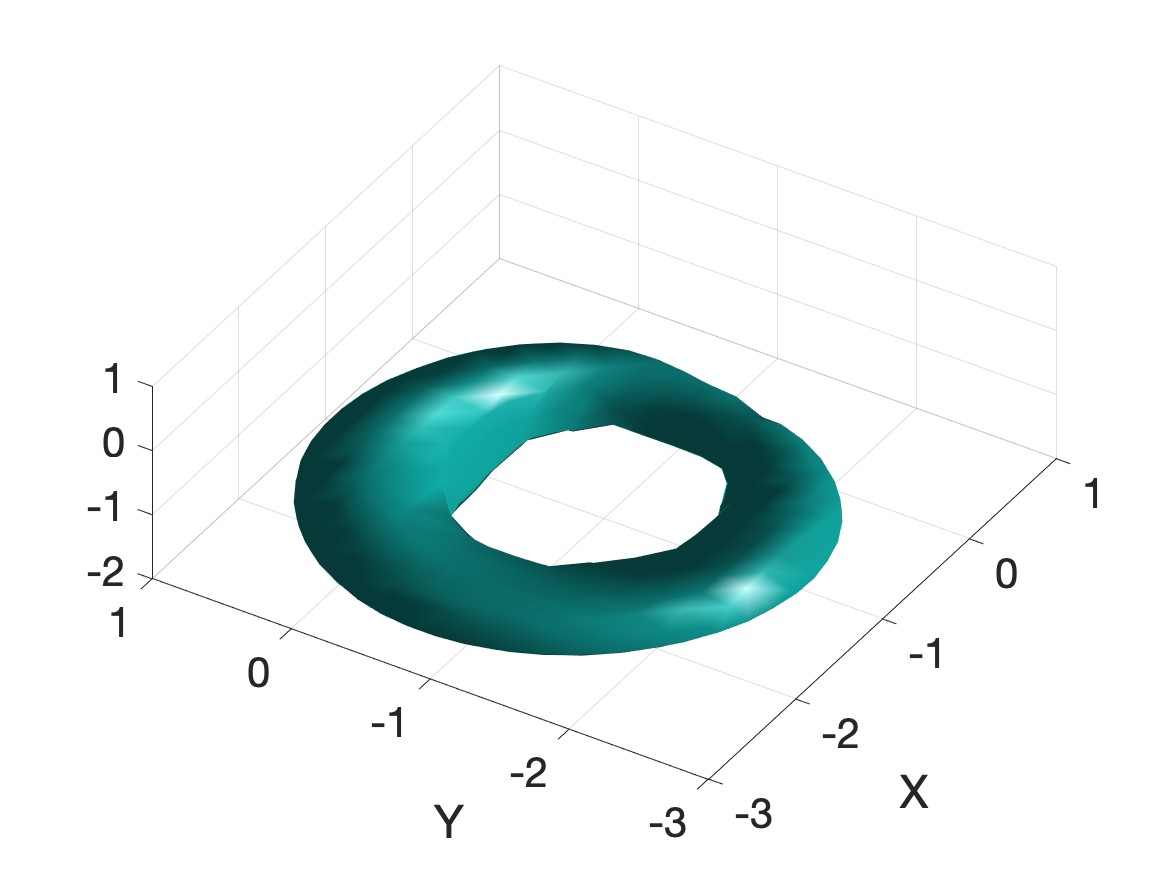}}
		
		\subfigure[The reconstructed letter `O' with multiple sources and a fixed frequency]{
		\includegraphics[width=0.35\textwidth]{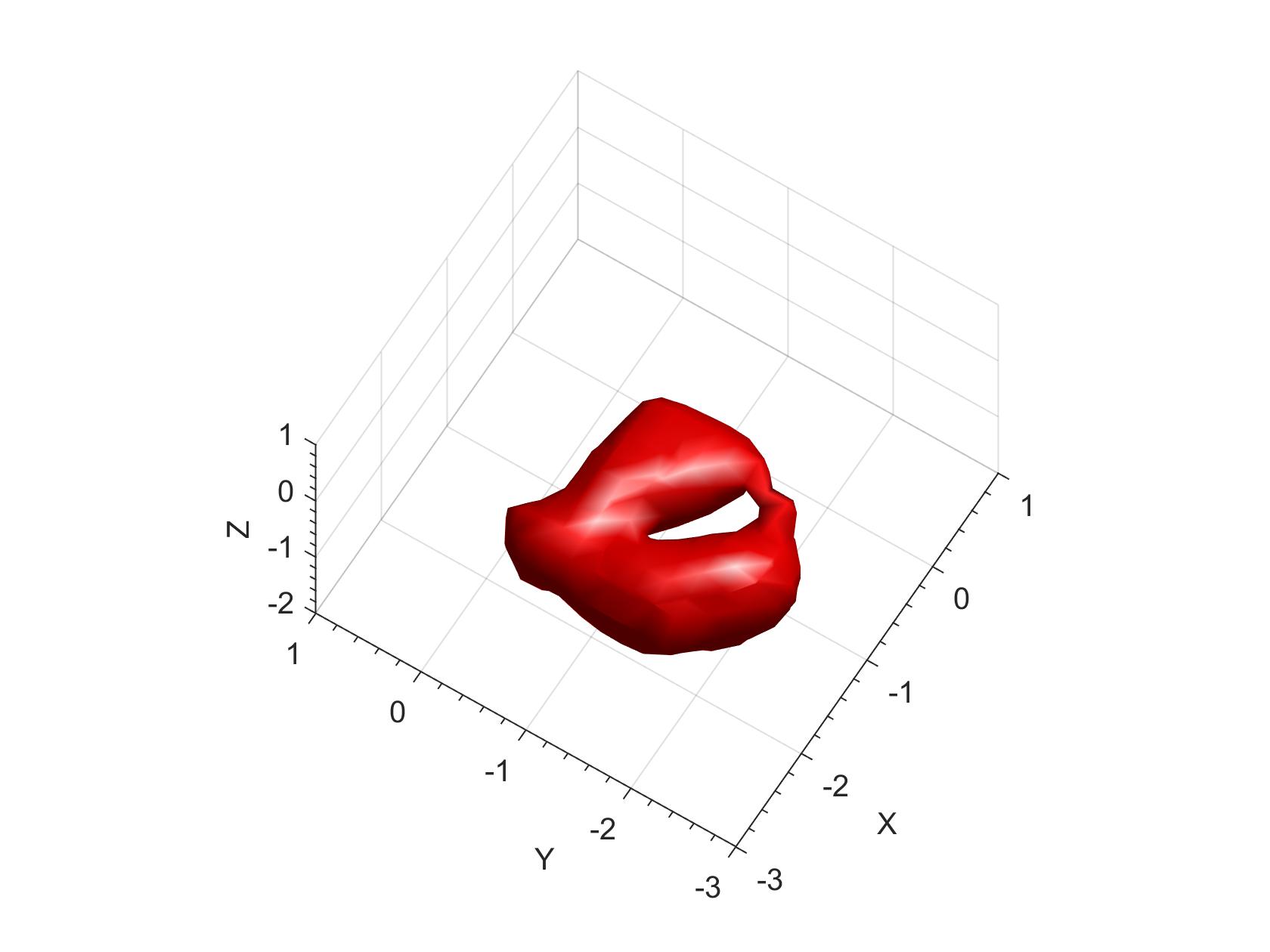}}
	\end{center}
	\caption{Example 5. The metallic letter `O'. (a) The real image of the O-shaped piece of metal. (b) The reconstructed dielectric constant function by the first method with a fixed source and multiple frequencies. (c) The reconstructed dielectric constant function by the second method with multiple sources and a fixed frequency. It is clear that the O shape with a void inside is produced better by the first method in comparison with the result of the second one.	\label{fig test 5}}

\end{figure}

Our numerical results are depicted in Figures \ref{fig test 1}--\ref{fig test 5} corresponding to five (5) examples that we have mentioned earlier in subsection \ref{subsec:setting}. In those figures, we present real photos of the experimental objects, and the reconstructed inclusions in three dimensions from two different computational approaches for comparison. Herewith, the first one is our current approach when using multi-frequency data and a fixed point source. The second approach is the one investigated in a series of publications  \cite{KhoaKlibanovLoc:SIAMImaging2020, VoKlibanovNguyen:IP2020, Khoaelal:IPSE2020} dealing with the context of multiple point sources and a fixed frequency.  

In all figures, we find that the first method reconstructs our inclusions with better shapes. Specifically, at this time the whole complicated shape of letters `U', `A' and `O' is visible when applying the first approach; see Figures \ref{fig test 3}, \ref{fig test 4} and \ref{fig test 5}. Also, in Figure \ref{fig test 2b}, the bottle of  water with its cap can be well interpreted, compared to Figure \ref{fig test 2c} in which the second approach is used. We, however, observe numerically that there is a main drawback of the current numerical approach. Our approach in this context (i.e. multiple frequencies and one source) does not give a decent value of dielectric constant. For instance, in the first test with the metallic cylinder, we report that the maximum value of the computed dielectric constant is only 1.0006, while by the second approach (i.e. multiple sources and one frequency), the obtained value is 18.72. Note that we, herewith, focus on the so-called ``appearing'' dielectric constant of metallic objects we have experimented numerically with in the previous publication \cite{Kuz2012}. In particular, the range of the appearing dielectric constant of metals is [10, 30].

\begin{table}
	\begin{centering}
		\begin{tabular}{|c|c|c|c|c|c|}
			\hline 
			Example & 1 & 2 & 3 & 4 & 5\tabularnewline
			\hline 
			Object & Metallic cylinder & Bottle of water & Wood ``U'' & Metallic ``A'' & Metallic ``O''\tabularnewline
			\hline 
			$c^{h}$ & 18.72 & 23.29 & 6.56 & 15.01 & 16.25\tabularnewline
			\hline 
			$c_{\text{true}}$ & 10 - 29 & 23.8 & 2 - 6 & 10 - 29 & 10 - 29\tabularnewline
			\hline 
			Reference & \cite{Kuz2012} & \cite{Thanh2015} & \cite{Clipper} & \cite{Kuz2012} & \cite{Kuz2012} \tabularnewline
			\hline 
		\end{tabular}
		\par\end{centering}
	\caption{Values of computed and true dielectric constants of examples 1--5. The values are taken from \cite{VoKlibanovNguyen:IP2020}. \label{tab:3}}
\end{table}

As one of important physical properties that one targets in landmine detection, shape of reconstructed inclusion is essential and can be helpful in classifying explosive devices in the battlefield. Our reconstruction results show that the perspective of multiple frequencies and one source being considered in this work does a good job to fulfill this property. It, indeed, produces a quite good shape of buried objects. With the same experimentally collected data used but different perspective (multiple sources and one frequency), the convexification method therein provides a high accuracy of computing the dielectric constant; see Table 2 in \cite{VoKlibanovNguyen:IP2020}. Note that having an accurate dielectric constant of the buried object is another essential physical property in landmine detection. Henceforth, when data set is allowed, it is our idea that one should combine these two perspectives to obtain decent reconstruction results in terms of both shape of the buried object and the dielectric constant. In the future work, we will find an appropriate approximate model for this interesting idea. In other words, a convexification method should be studied to come up with the perspective of multiple sources and frequencies.

\section{Summary} \label{sec:last}

In this paper, we have examined the numerical performance of our convexification method applied to a 3D coefficient inverse problem using experimental data. Our study focuses on imaging buried objects within a sandbox, simulating the detection of landmines on a battlefield.

Previously, we employed the convexification method with a setup involving multiple sources and a fixed frequency. This approach yielded highly accurate computations of the dielectric constants. Meanwhile, we observed that using multiple frequencies and a fixed source configuration improved the shape of the front surface of the experimental inclusions, a crucial physical property for detecting explosive devices.

Based on our current investigation, it is evident that combining these two configurations produces good reconstruction results in terms of both the shape of the buried object and the dielectric constant, provided the data set allows for such combination.

\begin{acknowledgments}
	V. A. Khoa thanks Dr. Darin Ragozzine (Brigham Young University, USA) for the recent support of his research career. L. H. Nguyen was supported by NSF Grant \#DMS-2208159.
\end{acknowledgments}

\section*{References}
\bibliographystyle{plain}
\bibliography{maikbibl}

\end{document}